\documentclass [12pt]{article}      
\usepackage{amsmath} 
\usepackage{latexsym}
\usepackage{pb-diagram}
\usepackage{amssymb}
\usepackage{bbm}
\usepackage{color}
\usepackage{graphicx}
\usepackage{inputenc}
\usepackage{hyperref}
\hypersetup{
    colorlinks,
    citecolor=blue,
    filecolor=blue,
    linkcolor=blue,
    urlcolor=blue
}
\hypersetup{linktocpage}
\usepackage{todonotes}
\usepackage{enumerate}

\newtheorem{theorem}{Theorem}
\newtheorem{corollary}[theorem]{Corollary}
\newtheorem{lemma}[theorem]{Lemma}
\newtheorem{definition}[theorem]{Definition}
\newtheorem{prop}[theorem]{Proposition}

\newtheorem{remark}[theorem]{Remark}

\newcommand{\rest}{\ensuremath{\upharpoonright}}

\newcommand{\mct}{{\ensuremath{\mathcal T}}}

\newcommand{\mcd}{{\ensuremath{\mathcal D}}}

\newcommand{\xbmt}{\ensuremath{(X,\mcb,\mu,T)}}
\newcommand{\ycns}{\ensuremath{(Y,\mcc,\nu,S)}}

\newcommand{\join}{\ensuremath{\bigvee}}

\newcommand{\qed}{{\nopagebreak \hfill $\dashv$ 
 \par\bigskip}}

\newcommand{\pf}{{\par\noindent{$\vdash$\ \ \ }}}

\newcommand{\poP}{\ensuremath{\mathbb P}}

\newcommand{\poR}{\ensuremath{\mathbb R}}

\newcommand{\la}{\langle}
\newcommand{\ra}{\rangle}

\newcommand{\mci}{\ensuremath{\mathcal I}}

\newcommand{\mcq}{\ensuremath{\mathcal Q}}
\newcommand{\mca}{\ensuremath{\mathcal A}}

\newcommand{\mcr}{\ensuremath{\mathcal R}}

\newcommand{\mcp}{\mathcal P}

	\newcommand{\poZ}{\mathbb Z}
	
	\newcommand{\mcw}{\mathcal W}
	\newcommand{\nn}{{\mathbb N}}

	\newcommand{\bfni}[1]{\noindent {{\bf{#1}}}}
	\newcommand{\mce}{{\mathcal E}}
	\newcommand{\inv}{{^{-1}}}

	\newcommand{\rev}[1]{\mathop{\rm rev}({#1})}

	\newcommand{\bk}{{\mathbb K}}
	
	\newcommand{\MPT}{\mbox{\bf MPT}}
	\newcommand{\pn}{{p_n}}
\newcommand{\pnpo}{{p_{n+1}}}
\newcommand{\qn}{{q_n}}
\newcommand{\qnpo}{{q_{n+1}}}

\newcommand{\kn}{{k_n}}

\newcommand{\mcc}{{\mathcal C}}

\newcommand{\zoo}{{[0,1)}}
\newcommand{\zoc}{{[0,1]}}

\newcommand{\mck}{{\mathcal K}}

\newcommand{\mcb}{\mathcal B}

\newcommand{\mcs}{{\mathcal S}}

\newcommand{\mudiff}{\mbox{Diff}^k(M, \mu)}
\newcommand{\sz}{\Sigma^\poZ}

\newcommand{\bt}{\mathbb T}
\newcommand{\hoo}[2]{{[{#1}/{#2}, ({#1}+1)/{#2}) }}
\newcommand{\rid}[1]{{\overline{\mcr}}_{#1}}

\newcommand{\mcu}{\mathcal U}

 \title{A symbolic representation for Anosov-Katok Systems}
 \author{Matthew Foreman\footnote{The first author would like to acknowledge partial support under  NSF award DMS 07010310}, Benjamin Weiss}

\begin{document}
 \maketitle

%
%

\begin{abstract}{This paper is the first of a series of papers culminating in the result that measure preserving diffeomorphisms of the disc or 2-torus are unclassifiable. It addresses another classical problem: which abstract measure preserving systems are realizable as smooth diffeomorphisms of a compact manifold? The main result gives symbolic representations of the Anosov-Katok diffeomorphisms. }
\end{abstract}
 \tableofcontents

 \section{Introduction}

      In 1932 J. von Neumann in \cite{vN} laid the foundations for ergodic theory. In it he expressed the likelihood that any abstract
      measure  preserving transformation (abbreviated to MPT in the paper)
        is isomorphic to a continuous MPT and
perhaps even to a differentiable one. Recall that two MPT's, $T$ and $S$,
are isomorphic if there is an invertible measure preserving  mapping between
the measure spaces which commutes with the actions of $T$ and $S$.
This brief remark eventually gave rise to one of the  outstanding problems in smooth dynamics, namely:

\begin{center}

 Does every
ergodic MPT with finite entropy have a smooth model?\footnote{In \cite{vN} on page 590,
``Vermutlich kann sogar zu jeder allgemeinen Str\"{o}mung eine isomorphe stetige Str\"{o}mung gefunden werden [footnote 13], vielleicht sogar eine stetig-differentiierbare, oder gar eine mechanische.
Footnote 13: Der Verfasser hofft, hierf\"{u}r demn\"{a}chst einen Beweis anzugeben."}

\end{center}

  By a smooth model it is meant an isomorphic copy of the MPT which is given by smooth diffeomorphism of a compact manifold preserving a measure equivalent to the volume element.
The  finite entropy restriction is required by a result of  A. G. Kushnirenko that showed that the entropy of any such  diffeomorphism must
be finite.
 An even more  basic problem which von Neumann  formulated in the same paper, was that of classifying
  all measure preserving transformations up to isomorphism.  This problem was solved long ago for several
   classes of  transformations that have special properties.  P. Halmos and von Neumann showed that ergodic MPT's with pure point spectrum
   are classified by the unitary equivalence of the associated unitary operators defined on the $L^2$ by the MPT, while A. N. Kolmogorov and D. Ornstein
   showed that Bernoulli shifts are classified by their entropy.

    One way to show  that not all finite entropy ergodic MPT's have a smooth model would be to show that their  classification is easier than the general
classification problem. Set theory provides a framework for a rigorous comparison of the complexity of different equivalence relations, and thus could potentially be a tool for settling this question.

Indeed, starting in the late 1990's a different type of result began to appear that used descriptive set theoretic techniques. These \emph{anti-classification} results  demonstrate in a rigorous way that positive classifications, such as those described above, are not possible.

The first is due to Beleznay and Foreman \cite{BF} who showed that the class of \emph{measure distal} transformations used in early ergodic theoretic proofs of Szemeredi's theorem is not a Borel set. Later Hjorth \cite{hj} introduced the notion of \emph{turbulence} and showed that there is no Borel way of attaching algebraic invariants to ergodic transformations that completely determine isomorphism. Foreman and Weiss \cite{FW} improved this results by showing that the conjugacy action of the measure preserving transformations is turbulent--hence no generic class  can have a complete set of algebraic invariants.

\smallskip

An ``anti-classification" theorem requires a precise definition of what a classification is. Informally a classification is a {method} of determining isomorphism between transformations perhaps by computing (in a liberal sense) other invariants for which equivalence is easy to determine.
The key words here are \emph{method} and \emph{computing}. For negative theorems, the more liberal a notion one takes for these words, the stronger the theorem. One natural notion is the Borel/non-Borel distinction. Saying a set $X$ or function $f$ is Borel is a loose way of saying that membership in $X$ or the computation of $f$ can be done using a countable (possibly transfinite) protocol whose basic input is membership in open sets. Say that $X$ or $f$ is \emph{not} Borel is saying that determining membership in $X$ or computing $f$ cannot be done with any countable amount of resources.

In the context of classification problems, saying that an equivalence relation $E$ on a space $X$ is \emph{not} Borel is
saying that there is no countable amount of information and no countable transfinite protocol for determining, for arbitrary
$x,y\in X$ whether $xEy$. {Any} such method must inherently use uncountable
resources.\footnote{Many well known classification theorems have as  immediate corollaries that the resulting equivalence
relation is Borel. An example of this is the Spectral Theorem, which has a consequence that the relation of Unitary
Conjugacy for normal operators is a Borel equivalence relation.}

In considering the isomorphism relation as a collection $\mathcal I$ of pairs $(S,T)$ of measure preserving transformations, Hjorth showed that $\mathcal I$ is not a Borel set. However the pairs of transformations he used to demonstrate this were inherently non-ergodic, leaving open the essential problem:

\begin{center}

Is isomorphism of ergodic measure preserving transformations Borel?

\end{center}

This question was answered in the negative by Foreman, Rudolph and Weiss in \cite{FRW}. This answer can be interpreted as saying that determining isomorphism between ergodic transformations is inaccessible to countable methods that use countable amounts of information.

   This  series of papers culminates in a result that--even restricted to the Lebesgue measure preserving diffeomorphisms of the 2-torus--the isomorphism relation is not Borel.

\begin{theorem}\label{big tuna} If $M$ is either the torus $\bt^2$, the disk $D$ or the annulus then the measure-isomorphism relation among pairs $(S,T)$ of measure preserving $C^\infty$-diffeomorphisms of $M$ is not a Borel set with respect to the $C^\infty$-topology.

\end{theorem}

\noindent{\bf What is in this paper?} The transformations built by Foreman, Rudolph and Weiss (\cite{FRW}) to prove the earlier result were based on odometers (in the sense that the Kronecker factor was an odometer). It is a well known open problem whether it is possible to have
 a smooth transformation on a compact manifold that has a non-trivial odometer factor. Thus proving the anti-classification theorem in the smooth context required constructing a different collection of hard-to-classify transformations and then showing that this collection could be realized smoothly.

 The new collection of transformations, the \emph{Circular Systems}, are defined as symbolic systems constructed using the \emph{Circular Operator}, a formal operation on words. This paper defines this class and then realizes them smoothly using the \emph{method of conjugacy} originating in a famous paper of Anosov and Katok.

 In fact, something much stronger is shown: Theorems \ref{verifying unique ergodicity} and \ref{circular as smooth} show that the Circular Systems exactly coincide with the isomorphism classes of the Anosov-Katok construction. We loosely summarize  \ref{verifying unique ergodicity} and \ref{circular as smooth}   as follows:

 \begin{theorem}[Main Result of this paper] \label{theorem 2}
 Let $T$ be an ergodic transformation on a standard measure space. Then the following are equivalent:
 \begin{enumerate}
        	 \item $T$ is isomorphic to an Anosov-Katok diffeomorphism.\footnote{Built using the untwisted method of conjugacy with some minor technical assumptions.}
	 \item $T$ is isomorphic to a (strongly uniform) circular system.
 \end{enumerate}
 \end{theorem}

This theorem shows that a broad class of transformations can be realized as Anosov-Katok diffeomorphisms. In fact we conjecture the following:

\begin{quotation}
\noindent {\bf Conjecture:} Suppose that $T$ is a zero-entropy ergodic transformation that has an irrational rotation of the circle as a factor. Then $T$ is isomorphic to a uniform circular system.
\end{quotation}
If the conjecture is true than all zero entropy ergodic transformations with a Liouvillean rotation factor can be realized as smooth transformations, and moreover, {every zero entropy ergodic transformation is a factor of an ergodic smooth transformation.}
\medskip

\bfni{Applications of Theorem \ref{theorem 2}:} Theorem \ref{theorem 2} is primarily useful in that it reduces questions about diffeomorphisms to combinatorial questions about symbolic shifts. Thus Theorem \ref{theorem 2} implies that any systematic way of building circular systems with given ergodic properties  automatically implies that there are ergodic measure preserving diffeomorphisms of the torus with those same properties.  

In \cite{global_structure} to this paper, the collection of circular shifts is endowed with a category structure and it is shown that this category is quite large. In particular it contains measure theoretically distal transformations of arbitrarily countable height.  By Theorem \ref{theorem 2} these are automatically realized as diffeomorphisms of the torus. 

Another sequel (\cite{FW}) reduces a complete analytic set (the \emph{ill-founded trees}) to the isomorphism relation for circular systems. By Theorem \ref{theorem 2} this automatically gives a reduction to the isomorphism problem for diffeomorphisms.

\bigskip

\bfni{Acknowledgements}

This work was inspired by the pioneering work of our co-author Dan Rudolph, who passed away before this portion of the grand project was undertaken. We owe an inestimable debt to J.P. Thouvenot who suggested using the Anosov-Katok method of conjugacy to produce our badly behaved transformations rather than attacking the ``odometer obstacle." The first author would like to thank Anton Gorodetsky for patiently explaining the original Anosov-Katok construction and its more contemporary versions as well as being a general reference on smooth dynamics.

\section{Preliminaries}\label{preliminaries}
In this section we establish some of the conventions we follow in this paper. There are many sources of background information on this including any standard text such as \cite{ergsurvey} or \cite{Peterson}.

\subsection{Measure Spaces}\label{abstract measure spaces} We will call separable non-atomic probability spaces \emph{standard measure spaces} and denote them $(X,\mcb, \mu)$ where $\mcb$ is the Boolean algebra of measurable subsets of $X$ and $\mu$ is a countably additive, non-atomic measure defined on $\mcb$. We will often identify two  members of $\mcb$ that differ
by a set of $\mu$-measure $0$ and seldom distinguish between $\mcb$ and the $\sigma$-algebra of classes of measurable sets modulo measure zero unless we are making a pointwise definition and need to claim it is well defined on equivalence classes.

\begin{remark}
Von Neumann proved that every standard measure space is isomorphic to $([0,1],\mcb,\lambda)$ where $\lambda$ is Lebesgue measure and $\mcb$ is the algebra of Lebesgue measurable sets.
\end{remark}

If $(X, \mcb, \mu)$ and $(Y, \mcc, \nu)$ are measure spaces, an isomorphism between $X$ and $Y$ is a bijection $\phi:X\to Y$ such that $\phi$ is measure preserving and both $\phi$ and $\phi^{-1}$ are measurable. We will ignore sets of measure zero when discussing isomorphisms; i.e.  we allow  the domain and range of $\phi$ to be  subsets of $X$  and $Y$ (resp.) of measure one.

A measure preserving system is an 4-tuple $\xbmt$ where $T:X\to X$ is a measure isomorphism. A \emph{factor map} 
between two measure preserving systems $\xbmt$ and $\ycns$ is a measurable, measure preserving  function 
$\phi:X\to Y$ such that $S\circ\phi=\phi\circ T$. A factor map is an isomorphism between systems iff $\phi$ is a 
measure isomorphism. As above we only require the domain and range of $\phi$ to have measure one, rather than that $\phi$ be one-to-one and onto.

\subsection{Partitions of measurable spaces}
We will  be concerned with ordered countable measurable partitions of measure spaces. 
An ordered countable measurable partition is a sequence $\mcp=\la P_n:n\in\nn\ra$ such that:
\begin{enumerate}
\item each $P_n\in \mcb$
\item if $n\ne m$ then $P_n\cap P_m=\emptyset$.
\item $\bigcup_nP_n$ has measure one.
\end{enumerate}
 We explicitly allow some of the $P_n$'s to be measure zero. The $P_n$'s will be called the \emph{atoms} of the partition.

If $\mcp=\la P_n:n\in \nn\ra$ and $\mcq=\la Q_n:n\in \nn\ra$ are two ordered partitions then the partition distance is defined as follows:
\[D_\mu(\mcp,\mcq)=\sum\mu(P_i\Delta Q_i).\]
We will frequently refer to  ordered countable measurable partitions simply as \emph{partitions}. A partition is finite iff for all large enough $n, \mu(P_n)=0$.  If we let $\poP_n$ be the space of partitions with $\le n$-atoms (i.e. for $m\ge n, \mu(P_m)=0$), then $(\poP_n, D_\mu)$ is a  connected space. 

If $\mcp$ and $\mcq$ are two partitions, then we say that $\mcq$ \emph{$\epsilon$-refines} $\mcp$ iff
the  atoms of $\mcq$ can be grouped into sets $\la S_n:n\in \nn\ra$ such that 
\[\sum_n\mu(P_n\Delta (\bigcup_{i\in S_n}Q_i))<\epsilon.\]
If $\mcp$ and $\mcq$ are partitions then $\mcq$ \emph{refines} $\mcp$ iff the atoms of $\mcq$ can be grouped into sets $\la S_n:n\in \nn\ra$ such that 
\[\sum_n\mu(P_n\Delta (\bigcup_{i\in S_n}Q_i))=0.\]
In this case we will write that $\mcq\ll\mcp$. A \emph{a decreasing sequence of partitions} is a sequence $\la \mcp_n:n\in\nn\ra$ such that for all $m<n, \mcp_n\ll \mcp_m$. If $A\in \mcb$ is a measurable set and $\mcp$ is a partition then we let $\mcp\rest A$  be the partition of $A$ defined as $\la P_n\cap A:n\in \nn\ra$.

\begin{definition}\label{definition of generation}
Let $(X,\mcb, \mu)$ be a measure space. We will say that a sequence of partitions $\la \mcp_n:n\in\nn\ra$ \emph{generates} (or generates $\mcb$) iff the smallest $\sigma$-algebra containing $\bigcup_n\mcp_n$ is $\mcb$ (modulo measure zero sets). If $T$ is a measure preserving transformation we will write $T\mcp$ for the partition $\la Ta:a\in \mcp\ra$. In the context of a measure preserving $T:X\to X$ we will say that a partition $\mcp$ is a \emph{generator} for $T$ iff $\la T^i\mcp:i\in \poZ\ra$ generates $\mcb$.

\end{definition}

We will be manipulating partitions of $\zoo$ and $\zoo\times \zoo$ in various ways so we develop some notation for 
doing so. We let \hypertarget{IQ}{$\mci_q$} be the partition of $[0,1)$ with atoms $\la \hoo{i}{q}:0\le i<q\ra$, and refer to $\hoo{i}{q}$ as $I^q_i$.\footnote{If $i>q$ then $I^q_i$ refers to $I^q_{i'}$ where $i'<q$ and $i'\equiv i \mod{q}$.} If $\mcp$ 
and $\mcq$ are partitions of spaces $X$ and $Y$ respectively, we let $\mcp\otimes \mcq$ be the partition of 
$X \times Y$ given by $\{P_i\times Q_j:i, j\in \nn\}$. To make this definition complete we need to fix in advance an arbitrary
ordering of $\nn\times \nn$ that is used to order of $\mcp\otimes \mcq$. Finally, we use the notation $I\otimes\mcq$ 
for the partition $\mcp\otimes \mcq$ where $\mcp$ has one element $I$.

If $T:X\to X$ and $\mcp=\la a_i:i\in I\ra$ is a partition of $X$ then the $(T,\mcp,n)$-\emph{name} of $s$ is $a_i$ if and only if $T^n(x)\in a_i$. If $T$ is invertible then the $(T,\mcp)$-name is $s\in \mcp^\poZ$ if and only if for all $n\in \poZ, T^n(x)\in s(n)$. We suppress $\mcp$ and/or $T$ if either is obvious from context.


\section{Presentations of Measure Preserving Systems}\label{presentations}
Measure preserving systems occur naturally in many guises with diverse topologies. As far as is known, the Borel/non-Borel distinction for dynamical properties is the same in each of these presentations and many of the presentations have the same generic classes. (See the forthcoming paper \cite{models}.)

In this section we briefly review the properties of the presentations relevant to this paper. These are: 
abstract invertible preserving systems, smooth transformations preserving volume elements and symbolic 
systems. 

\subsection{Abstract Measure Preserving systems}
As noted in section \ref{abstract measure spaces} every standard measure space is isomorphic to the unit interval with Lebesgue measure. Hence every invertible measure preserving transformation of a standard measure space is isomorphic to an invertible Lebesgue measure preserving transformation on the unit interval. 

In accordance with the conventions of \cite{ergsurvey} we denote the collection of measure preserving transformations of $\zoo$ by \MPT.\footnote{Recently several authors have adopted the notation $Aut(\mu)$ for the same space.} We note that two measure preserving transformations are identified if they are equal on sets of full measure. 

We can associate to each invertible measure preserving transformation $T\in \MPT$ a unitary operator $U_T:L^2([0,1])\to L^2([0,1])$ by defining $U(f)=f\circ T$.  In this way \MPT\ can be identified with a closed subgroup of the unitary operators on $L^2([0,1])$ with respect to the weak operator 
topology
on the space of  unitary 
transformations. This makes \MPT\ into a Polish space. We will call this the \emph{weak topology} on \MPT\ (See \cite{Halmos}).

A concrete description of the topology can be given as follows: Let $S\in \MPT$, $\mcp$ be a finite measurable partition and $\epsilon>0$. Define 
\begin{eqnarray*}
N(S, \mcp, \epsilon)=\{T\in \MPT: \sum_{A\in \mcp}\lambda(TA\Delta SA)<\epsilon\}
\end{eqnarray*}
If $\la \mcp_n:n\in\nn\ra$ is a generating sequence of partitions for $\mcb$, then 
 $\{N(S, \mcp_n, \epsilon):S\in \MPT, n\in\nn, \epsilon>0\}$ generates the weak operator topology on \MPT.
 
 We will denote the ergodic transformations belonging to $\MPT$ by $\mce$. Halmos (\cite{Halmos}) showed that 
 $\mce$ is a dense $\mathcal G_\delta$ set in $\MPT$. In particular the weak topology makes $\mce$ into a Polish subspace of $\MPT$.
 
 The following is easy to check:
 \begin{lemma}\label{criterion for weak convergence}
 Let $\la T_n:n\in\nn\ra$ be a sequence of measure preserving transformations and $\la \mcp_n:n\in\nn\ra$ be a generating sequence of partitions. Then the following are equivalent:
 \begin{enumerate}
  \item The sequence $\la T_n:n\in \nn\ra$ converges to an invertible measure preserving system in the weak topology.
  \item For all measurable sets $A$, for all $\epsilon>0$ there is an $N$ for all $n,m>N$ and $i=\pm 1$ we have 
  $\mu(T^i_nA\Delta T^i_mA)<\epsilon$. 
     \item For all $\epsilon>0, p\in\nn$ there is an $N$ for all $m,n>N$ 
     \[\sum_{A\in \mcp_p, i=\pm 1}\mu(T^i_nA\Delta T^i_mA)<\epsilon\]
\end{enumerate}
 In case the sequence $\la T_n:n\in\nn\ra$ converges then we can identify the limit as the unique $T$ such that for all measurable sets $A$, 
\[\mu(T_nA\Delta TA)\to 0.\]
 \end{lemma}
  
 There is another topology on the collection of measure preserving transformations of $X$ to $Y$  for measure spaces $X$ and $Y$. If $S, T:X\to Y$ are measure preserving transformations, the \emph{uniform distance} between $S$ and $T$ is defined to be:
\[d_U(S, T)=\mu\{x:Sx\ne Tx\}.\]
This topology refines the weak topology and is a complete, but  not a separable topology.

\subsection{Diffeomorphisms}
Let $M$ be a $C^k$-smooth compact finite dimensional manifold and $\mu$ be a standard measure on $M$ determined by a smooth volume element. 
For each $k$ there is a Polish topology on the $k$-times differentiable homeomorphisms of $M$, the $C^k$-topology. The $C^\infty$-topology is the coarsest topology refining the $C^k$-topology for each $k\in \nn$. 
It is also a Polish topology and a sequence of $C^\infty$-diffeomorphisms converges in the $C^\infty$-topology if and only if it converges in the $C^k$-topology for each $k\in \nn$. The $C^\infty$ topology is also a Polish topology and we will sometimes use a Polish metric $d^\infty$ on the diffeomorphisms inducing this topology.

The collection of $\mu$-preserving  diffeomorphisms forms a closed nowhere dense set in the $C^k$-topology on the $C^k$-diffeomorphisms, and as such inherits a Polish topology
\footnote{One can also consider the space of measure preserving homeomorphisms with the $\|\ \|_\infty$ topology, which behaves in some ways similarly.}.  
We will denote this space by $\mudiff$.

The measure preserving diffeomorphisms of a compact manifold can also be endowed with the weak topology, which is coarser than the $C^k$-topology. To see that the weak topology is coarser than the $C^k$-topologies, note that if $M$ is compact and has dimension $n$, then $M$ has a countable generating sequence  of finite partitions into ``half-open" sets whose boundaries are finite unions of submanifolds of dimension less than  $n$. Let $\mcp$ be such a partition. Then the boundaries of the elements of $\mcp$ all have measure zero and if $S$ and $T$ are close in the $C^k$-topology, then $S$ and $T$ take the boundaries to very similar places. In particular, $S\mcp$ and $T\mcp$ don't differ very much.

One can also consider the space of abstract $\mu$-preserving transformations on $M$ with the weak topology. In 
\cite{bjork_thesis} it is shown that the collection of a.e.-equivalence classes of smooth transformations form  a 
$\Pi^0_3$-set (${\mathcal G}_{\delta\sigma\delta}$) in \MPT(M), and hence the collection has the  Property of Baire. In particular, by invariance it is either meager or comeager.

\subsection{Symbolic Systems}
\label{symbolic shifts}
Let $\Sigma$ be a countable or finite alphabet endowed with the discrete topology. Then $\Sigma^\poZ$ can be given the product topology, which makes it into a separable, totally disconnected space that is compact if $\Sigma$ is finite.  
\medskip

\bfni{Notation:} If $u=\la \sigma_0, \dots \sigma_{n-1}\ra\in \Sigma^{<\infty}$ is a finite sequence of elements of $\Sigma$, then we denote the cylinder set based at $k$ in $\Sigma^\poZ$ by writing $\la u\ra_k$. If $k=0$ we abbreviate this and write $\la u\ra$. 
Explicitly: $\la u\ra_k=\{f\in \Sigma^\poZ: f\rest[k,k+n)=u\}$. The collection of cylinder sets form a base for the product topology on $\Sigma^\poZ$, thus we frequently refer to them as ``basic open sets."
\medskip

\noindent The shift map:
\[sh:\Sigma^\poZ\to \Sigma^\poZ\]
defined by setting $sh(f)(n)=f(n+1)$ is a homeomorphism. If $\mu$ is a shift-invariant Borel measure then the resulting measure preserving system $(\Sigma^\poZ, \mcb,\mu, sh)$ is called a \emph{symbolic system}. The closed support of $\mu$ is a shift-invariant closed subset of $\Sigma^\poZ$ called a \emph{symbolic shift} or \emph{sub-shift}. 

We can construct symbolic shifts from arbitrary measure preserving systems as follows: If $\xbmt$ is a measure preserving system and $\mcp=\{A_i:i\in I\}$ is a measurable partition (where $I$ is countable or finite). Let $\Sigma=\{a_i:i\in I\}$. then we can define a map
\begin{eqnarray*} \phi:X\to \Sigma^\poZ
\end{eqnarray*}
by setting $\phi(x)(n)=a_i$ iff $T^nx\in A_i$. 

The map $\phi$ induces an invariant Borel measure $\nu=\phi^*\mu$ on $\Sigma^\poZ$ by setting $\nu(B)=\mu(\phi^{-1}(B))$. The resulting invariant measure makes $(\Sigma^\poZ, \mcc, \nu, sh)$ into a factor of $\xbmt$ with factor map $\phi$. Since $X$ is standard, if $\mathcal P$ generates then $\phi$ is an isomorphism.  

\begin{remark}\label{noise}
{We will use the fact that we can systematically change symbols in some positions of letters in   $x\in\Sigma^\poZ$ to get a new element $x'\in \Sigma^\poZ$ as long as the change is equivariant with the shift and the map $x\mapsto x'$ is one to one. Because the change is one to one we can copy over the measure $\nu$ to a measure $\nu'$ so that the  resulting measure on $(\Sigma)^\poZ$ will define an isomorphic system.}
\end{remark}

{ \bfni{Notation:} For a word $w\in \Sigma^{<\nn}$ we will write  $|w|$ for the length of $w$.}
\smallskip

We want to be able to unambiguously parse elements words and elements of symbolic shifts. For this we will use construction sequences consisting of uniquely readable words.
\begin{definition}
Let $\Sigma$ be a alphabet and $\mcw$ be a collection of finite words in $\Sigma$. Then $\mcw$ is \emph{uniquely readable} iff whenever $u, v, w\in \mcw$ and $uv=pws$ then either $p$ or $s$ is the empty word.
\end{definition}

Symbolic shifts are often described intrinsically by giving a collection of words that constitute a clopen base for the support of an invariant measure.  
Fix a alphabet $\Sigma$,  and a sequence of uniquely readable collections of    words $\la\mcw_n:n\in\nn\ra$ with the properties that:
\begin{enumerate}
\item for each $n$ all of the words in $\mcw_n$ have the same length $q_n$,
\item each $w'\in \mcw_{n+1}$ contains each $w\in\mcw_n$ as a subword,
\item there is a summable sequence $\la \epsilon_n:n\in\nn\ra$ of positive numbers such that for each $n$, every word $w\in \mcw_{n+1}$ can be uniquely parsed into segments 
\[u_0w_0u_1w_1\dots  w_lu_{l+1}\]
 such that each $w_i\in \mcw_n$, $u_i\in \Sigma^{<\nn}$ and for this parsing
\begin{equation}\label{small boundary numeric} {\sum_i|u_i|\over q_{n+1}}<\epsilon_{n+1}.
\end{equation}
\end{enumerate}
\begin{definition}A sequence $\la \mcw_n:n\in\nn\ra$ satisfying items 1.)-3.) will be called a \emph{construction sequence}.
\end{definition}
 Define $\bk$ to be the collection of $x\in \Sigma^\poZ$ such that every finite contiguous subword of $x$ occurs inside some $w\in \mcw_n$. Then $\bk$ is a closed shift-invariant subset of $\sz$ that is compact if $\Sigma$ is finite. 
 \medskip

\begin{definition}\label{definition of uniform}
Let $\la\mcw_n:n\in\nn\ra$ be a construction sequence.  Then $\la \mcw_n:n\in\nn\ra$ is  \hypertarget{uniform}{\emph{uniform}} 
if 
there are $\la d_n:n\in\nn\ra$, where $d_n:\mcw_n\to (0,1)$ and a sequence $\la \epsilon_n:n\in\nn\ra$ going to zero such that for 
each $n$ all  words $w\in \mcw_n$ and $w'\in \mcw_{n+1}$ if $f(w,w')$ is the number of $i$ such that  $w=w_i$
\begin{equation}\label{uniform occurrance}
\sum_{w\in \mcw_n}\left|{f(w,w')\over q_{n+1}/q_n}-d_n(w)\right|<{\epsilon_{n+1}}.
\end{equation}
\end{definition}
\medskip

The words $u_i$ are often called \emph{spacers}. The $d_n$ are target values for the densities of 
$n$-words in $n+1$ words. The \emph{uniformity} is that each $n$-word occurs nearly the same number of times in every $n+1$-word. If 
$\bk$ is built from a uniform construction sequence we will call $\bk$ a \emph{uniform symbolic system}. 

If  $f(w,w')$ is a constant (depending on $n$) for all $w\in \mcw_n,w'\in \mcw_{n+1}$  we can take $d_n(w)={f(w,w')\over q_{n+1}/q_n}$ and satisfy definition \ref{definition of uniform}.  In this case  we call the construction sequence and  $\bk$ \emph{strongly uniform}.

 If $\bk\subset \Sigma^\poZ$ is a symbolic system, then an element $x\in \bk$ is a function $x:\poZ\to \Sigma$. If $I$ is a finite or infinite interval in $\poZ$, then we write $x\rest I$ for the function $x$ restricted to this interval. In our constructions we will  restrict our measures to a natural set:

\begin{definition}\label{def of S} Suppose that $\la \mcw_n:n\in\nn\ra$ is a construction sequence  for a symbolic system $\bk$ with each $\mcw_n$ uniquely readable. Let
\hypertarget{def of s}{$S$} be the collection $x\in \bk$ such that there are  sequences of natural numbers 
$\la a_m: m\in\nn\ra$, $\la b_m: m\in\nn\ra$ going to infinity such that for all large enough $m$,
$x\rest [-a_m, b_m)\in \mcw_m$.
\end{definition}


\noindent Note that  $S$ is a dense shift-invariant $\mathcal G_\delta$ set.

\begin{lemma} \label{unique ergodicity}
Fix a construction sequence $\la\mcw_n:n\in\nn\ra$ for a symbolic system $\bk$ in a finite alphabet $\Sigma$. Then:
\begin{enumerate}
\item $\bk$ is the smallest shift-invariant closed subset of $\Sigma^\poZ$ such that for all $n$,  and $w\in\mcw_n$, $\bk$ has non-empty intersection with the basic open interval $\la w\ra\subset \Sigma^\poZ$.
\item \label{stereotype}
Suppose that $\la \mcw_n:n\in\nn\ra$ is a uniform construction sequence. 
Then there is a unique non-atomic shift-invariant measure  $\nu$ on $\bk$ concentrating on $S$ and this $\nu$ is ergodic.

\end{enumerate}
\end{lemma}
\pf Item 1 is clear from the definitions.
To see item 2, fix a measure $\nu$ concentrating on $S$. It suffices to show that the $\nu$-measures of sets of the form $\la u\ra_0$ for $u\in \mcw_k$ are uniquely determined by $\la \mcw_n:n\in\nn\ra$.
Fix a $u\in \mcw_k$ for some $k$. 
By the Ergodic Theorem it suffices to show that for all $\epsilon>0$ and all large enough 
$n$ if $w', w''\in \mcw_{n+1}$, then the proportion of occurrences of $u$ among the $k$-words in $w'$ is within $\epsilon$ of the proportion of the $k$-words occurring in $w''$.

For each $w\in \mcw_n$, let $\lambda_w(u)$ be the proportion of occurrences of $u$ among the $k$-words occurring in $w$. Then the proportion of occurrences of $u$ among the $k$-words in $w'$ is approximated up to the proportion of $w'$ taken up by spacers (which is summably small) by
\[\sum_{w\in \mcw_n}\lambda_w(u){f_n(w,w')\over{q_{n+1}/q_n}}\]
and the similar approximation holds for $w''$. Computing:
\begin{eqnarray*}
\left|\sum_{w\in\mcw_n}\lambda_w(u){f_n(w,w')\over{q_{n+1}/q_n}}-\sum_{w\in\mcw_n}\lambda_w(u){f_n(w,w'')\over{q_{n+1}/q_n}}\right|&\le&\\
\sum_{w\in\mcw_n}\lambda_w(u)\left|{f_n(w,w')\over{q_{n+1}/q_n}}-{f_n(w,w'')\over{q_{n+1}/q_n}}\right|&\le&\\
\sum_{w\in\mcw_n}\left|{f_n(w,w')\over{q_{n+1}/q_n}}-{f_n(w,w'')\over{q_{n+1}/q_n}}\right|&\le&\\
\sum_{w\in\mcw_n}\left(\left|{f_n(w,w')\over{q_{n+1}/q_n}}-d_n(w)\right|+\left|d_n(w)-{f_n(w,w'')\over{q_{n+1}/q_n}}\right|\right)&\le&\\
&\le&2\epsilon_{n+1}.
\end{eqnarray*}
Taking $n$ large enough we have shown that $\nu(\la u\ra_0)$ is uniquely determined.  Since there is a unique measure on $S$, that measure must be ergodic.
%
\qed

\begin{remark}We make two remarks about Lemma \ref{unique ergodicity}.

\begin{enumerate}
\item If $X$ is a Polish space, $T:X\to X$ is a Borel automorphism and  $D$ is a $T$-invariant Borel set with a unique $T$-invariant measure on $D$, then that measure must be ergodic.

\item If $(\bk, sh)$ is an arbitrary symbolic shift then its inverse is $(\bk, sh^{-1})$, where $sh^{-1}(f)(n)=f(n-1)$.  If $x$ is in $\bk$, we define the 
reverse  of $x$ by setting $\rev{x}(k)=x(-k)$. We can view $(\bk, sh^{-1})$ as the symbolic system $(\rev{\bk}, sh)$, where 
$\rev{\bk}$ consists of all of the reverses  of elements of $\bk$.

\item {Assuming the hypothesis of Lemma \ref{unique ergodicity},
 the proof also shows that there is a unique non-atomic shift-invariant measure on $\rev{S}$ and that for this measure, which we denote $\nu\inv$, we have $\nu(\la w\ra)=\nu\inv(\la \rev{w}\ra)$.}

 \end{enumerate}
\end{remark}


 
\section{Circular Symbolic Systems}\label{odometer based and circular symbolic systems}

We now define a class of symbolic shifts that we call \emph{circular systems}. The main result of this paper is that the circular systems are the symbolic representations of the smooth diffeomorphisms defined by the Anosov-Katok method of conjugacies. The construction sequences of circular systems have quite specific combinatorial properties that will be important in the sequels to our understanding of the Anosov-Katok systems and their centralizers.

These symbolic systems are built from construction sequences $\la \mcw_n:n\in\nn\ra$ where $\mcw_{n+1}$ is the result of applying an abstract operation $\mcc$ to sequences of words from $\mcw_n$. We call these systems \emph{circular}  because they are closely tied to the behavior of rotations by a convergent sequence of rationals $\alpha_n=p_n/q_n$.  The rational rotation by $p/q$ permutes the $1/q$ intervals of the circle cyclically along a sequence determined by some numbers $j_i=_{def}p^{-1}i$ (mod $q$).\footnote{We assume that $p$ and $q$ are relatively prime and the exponent $-1$ denotes the multiplicative inverse of $p$ mod $q$.} 
To have a symbolic representation of an Anosov-Katok diffeomorphism, one must be able to describe how the intervals of length $1/q_{n+1}$ are permuted by  addition of $\alpha_{n+1}$ mod(1) in terms of the intervals of length $1/q_n$. The abstract symbolic operation $\mcc$ does this. We explain this in detail in section
\ref{symbolic representations of AK}.
 Theorem \ref{circular as smooth}, which says that all strongly uniform systems built with $\mcc$ and a suitable coefficient sequence can be realized as measure preserving diffeomorphisms, is the  part of this paper that we use in the sequels to construct measure preserving diffeomorphisms  with complicated combinatorial structures. Theorem \ref{verifying unique ergodicity} is the companion to Theorem \ref{circular as smooth}. It says that all (unwisted, strongly uniform) Anosov-Katok diffeomorphisms can be represented by circular symbolic systems.

Let $k, l, p, q$ be positive natural numbers with $p$ and $q$ relatively prime. Set $j_i\equiv_q(p)^{-1}i$ with $j_i<q$. 
It is easy to verify that:
\begin{equation}\label{reverse numerology} 
q-j_i=j_{q-i} 
\end{equation}

The operation  $\mcc$ is defined on sequences  $w_0, \dots w_{k-1}$ of words in a alphabet $\Sigma\cup \{b, e\}$ (where we assume that neither $b$ nor $e$ belong to $\Sigma$) by setting:

\begin{equation}\mcc(w_0,w_1,w_2,\dots w_{k-1})=\prod_{i=0}^{q-1}\prod_{j=0}^{k-1}(b^{q-j_i}w_j^{l-1}e^{j_i}). \label{definition of C}
\end{equation} 

\begin{remark}\label{word length}
We remark:
\begin{itemize}
\item Suppose that each $w_i$ has length $q$, then the length of  $\mcc(w_0, w_1, \dots w_{k-1})$ is $klq^2$. 
\item If $e$ occurs in $\mcc(w_0, \dots w_{k-1})$ then there is an occurrence of $b$ to the left of it.
\item {Suppose that $n<m$ and $b$ occurs at position $n$ and $e$ occurs at position  $m$ and neither occurrence is in a $w_i$. Then there must be some $w_i$ occurring between $n$ and $m$.}
\end{itemize}
\end{remark}

The following unique readability lemma is used to show that many construction sequences for circular systems are strongly uniform. We will also use it when we  show  that our symbolic representations of diffeomorphisms come from generating partitions.

\begin{lemma}\label{unique readability}
Suppose that $\Sigma$ is a finite or countable alphabet and that  $u_0, \dots $ $u_{k-1}$, $v_0, \dots v_{k-1}$ and $w_0 \dots w_{k-1},$ are words in the alphabet $\Sigma\cup \{b, e\}$ of some fixed length $q<l/2$. Let 
\begin{eqnarray*}u&=&\mcc(u_0, u_1, \dots u_{k-1}) \\
v&=&\mcc(v_0, v_1, \dots v_{k-1})\\
w&=&\mcc(w_0, w_1, \dots w_{k-1}).
\end{eqnarray*} Suppose that $uv$ is written as $pws$ where $p$ and $s$ are words in $\Sigma\cup \{b, e\}$. Then either $p$ is the empty word and $u=w, v=s$ or $s$ is the empty word and $u=p, v=w$. 
\end{lemma}
\pf We note that the map $i\mapsto j_i$ is one-to-one. Hence each location in the word of length $klq^2$ is uniquely determined by the lengths of nearby sequences of $b$'s and $e$'s.\qed

In fact something stronger is true: if $\sigma\in \Sigma$ occurs at place $m$ in $w$ then $m$ is uniquely determined by $w_0, w_1, \dots w_{k-1}$ and the $q^l/2 +1$ letters on either side of $\sigma$. 

We now describe how to use the $\mcc$ operation to build a collection of symbolic shifts. Our systems will be defined using a sequence of natural number parameters $ k_n$ and $l_n$ that are fundamental to the version of the Anosov-Katok construction as presented in \cite{katoks_book}.

 The numbers {$\la l_n\ra$} will be assumed to go to infinity quite rapidly.\footnote{In particular in what follows we will assume that  $\sum{1\over l_n}$ is finite.} From the $k_n$ and $l_n$ we define other sequences of numbers: $\la p_n, q_n, \alpha_n:n\in\nn\ra$ (with more defined later). 
We begin by letting $p_0=0$ and $q_0=1$  and inductively set 
\begin{eqnarray}\label{pn and qn}
\pnpo=\pn\qn\kn l_n +1 & & \qnpo=\kn l_n\qn^2.
\end{eqnarray}
Thus $p_1=1$ and $q_1=k_0l_0$. Letting $\alpha_n=p_n/q_n$ we see that:
\[{p_{n+1}\over q_{n+1}}=\alpha_n+{1\over k_nl_nq_n^2}.\]

We note that $p_n$ and $q_n$ are relatively prime for $n\ge 1$ and hence it makes sense to define an integer $j_i$ with $0\le j_i<q_n$ by setting:\footnote{{For $q_0=1$, $\poZ/q_0\poZ$ has one element, $[0]$, so $p_0\inv=p_0=0$.} Also, formally $j_i$ should have a notation indicating that it depends on $n$, i.e. $j_i^n$. We neglect this to reduce notational complexity. }
\begin{eqnarray}
j_i=(p_n)^{-1}i \mod q_n.  \label{j sub i}
\end{eqnarray}

 Let $\Sigma$ be a non-empty finite or countable alphabet. We will construct the systems we study by building collections of words $\mcw_n$ in $\Sigma\cup \{b, e\}$ by induction. We set $\mcw_0=\Sigma$.  Having built $\mcw_n$ we choose a set of \emph{prewords} $P_{n+1}\subseteq (\mcw_{n})^{k_n}$ and  form $\mcw_{n+1}$ by taking all words of the form $\mcc(w_0,w_1\dots w_{k_n-1})$  with $(w_0, \dots w_{k_n-1})\in P_{n+1}$.\footnote{Passing from $\mcw_n$ to $\mcw_{n+1}$ we use $\mcc$ with parameters $k=k_n, l=l_n, p=p_n$ and $q=q_n$. The $j_i$ is $(p_n)^{-1}i$ modulo $q_n$. Strictly speaking we should probably write $\mcc_n$ for the operation $\mcc$ at stage $n$ that uses these parameters and write $j_i$ as $j_i^n$. By Remark \ref{word length}, the length of each of the words in $\mcw_{n+1}$ is $q_{n+1}$.} It follows from Lemma \ref{unique readability} that each $\mcw_n$ is uniquely readable.
 \medskip
 
  \hypertarget{strong unique readability}{\bfni{Strong Unique Readability Assumption:}}
  Let $n\in \nn$, and view $\mcw_n$ as a collection $\Lambda_n$ of letters. We will say that $P_{n+1}$ satisfies strong unique readability if and only if when  viewing each element of $P_{n+1}$  as a word  with letters in   $\Lambda_n$, the resulting collection of $\Lambda_n$--words is uniquely readable.
  
  \begin{definition}\label{ccs} A construction sequence $\la \mcw_n:n\in\nn\ra$  will be called \emph{circular} if 
  it is built using the operation $\mcc$, a circular coefficient sequence and each 
  $P_{n+1}$ satisfies the strong unique readability assumption.   \end{definition}

%
We now show that strongly uniform circular systems are uniform.
\begin{lemma}\label{uniformity of circular systems}
Suppose  $\la \mcw_n:n\in\nn\ra$ is a circular construction sequence such that:
\begin{enumerate}  \item $\sum 1/l_n$ is finite and
 \item for each $n$ there is a number $f_n$ such that each word $w\in \mcw_n$ occurs exactly $f_n$ times in each word in $P_{n+1}$.

\end{enumerate}Then $\la \mcw_n:n\in\nn\ra$ is strongly uniform.
\end{lemma}

\pf 
%

For each $w\in \mcw_n, w'\in \mcw_{n+1}$, if we set 
$f(w,w')=f_nq_n(l_n-1)$ and $d_n=\left({f_n\over k_n}\right)\left(1-\left({1\over l_n}\right)\right)$ then:
\begin{eqnarray*}
{f(w,w')\over{q_{n+1}/q_n}}&=&f_nq_n(l_n-1)\left({q_n\over k_nl_nq_n^2}\right)\\
&=&\left({f_n\over k_n}\right)\left(1-\left({1\over l_n}\right)\right)\\
&=&d_n.
\end{eqnarray*}
\qed

\begin{definition}\label{circular definition}
A symbolic shift $\bk$ constructed from a  circular construction sequence  will be called a \emph{circular system}. If $\bk$ is constructed from a (strongly) uniform circular construction sequence then we will say that $\bk$ is a \emph{(strongly) uniform circular system}.
\end{definition}

 Lemma \ref{unique ergodicity} gives a characterization of the support of a uniform circular system and shows that there is a unique shift-invariant measure on the set $S$.

\begin{definition}
Suppose that $w=\mcc(w_0,w_1,\dots w_{k-1})$.  Then $w$ consists of blocks of $w_i$ repeated $l-1$ times, together with some $b$'s and $e$'s that are not in the $w_i$'s. The \emph{interior} of $w$ is the portion of $w$ in the $w_i$'s. {The remainder of  $w$ consists of blocks of the form $b^{q-j_i}$ and $e^{j_i}$. We call this portion the  \emph{boundary} of $w$.} 

In a block of the form $w_j^{l-1}$ the first and last occurrences  of $w_j$ will be called the  \emph{boundary portion of the block} $w_j^{l-1}$. The other occurrences will be the \emph{interior} occurrences. 
\end{definition}
{We note that the boundary consists of sections of $w$ made up of $b$'s and $e$'s. However not all $b$'s and 
$e$'s occurring in $w$ are in the boundary, as they may be part of a power $w_i^{l_i-1}$.}

 The boundary of $w$ constitutes a small portion of the word:

\begin{lemma}\label{stabilization of names 1} The proportion of the word $w$ written in
equation \ref{definition of C} that belongs to its boundary is $1/l$. Moreover the proportion of the word that is within $q$ letters of boundary  is $3/l$. 
\end{lemma}

We now characterize the set $S\subset \bk$  for circular systems and show a strong unique ergodicity result.

\begin{lemma}\label{dealing with S} Let $\bk$ be a circular system. Then 
\begin{enumerate}
\item Let $\nu$ be a shift-invariant measure on $\bk$, then  $\nu$ concentrates on $S$ iff  $\nu$ concentrates on  the collection of $s\in \bk$ such that $\{i:s(i)\notin \{b, e\}\}$ is 
unbounded in both $\poZ^-$ and $\poZ^+$.
\item Suppose that $\bk$ is a uniform circular system. If $\nu$ is a non-atomic shift-invariant measure on $\bk$ then $\nu(S)=1$.  In particular, there is a unique non-atomic shift-invariant measure on $\bk$. 
\end{enumerate}

\end{lemma}

\pf
To start the proof we note that if $e$ occurs at $m_0$ in $s\in S$, then there is an $m_1<m_0$ and a $\sigma\in \Sigma$ such that $\sigma$ occurs at $m_1$ in $s$. Similarly if $b$ occurs at $m_0$ there is an $m_1>m_0$ and a $\sigma\in \Sigma$ such that $\sigma$ occurs at $m_1$. To see this fix an occurrence $m_0$ of $e$ in $s$. (The argument for $b$ is symmetric.) Let $n$ be the smallest natural number such that there is a $w\in \mcw_{n+1}$ occurring in $s$ on the interval $[a,b]$ with $a\le m_0\le b$. Then $w=\mcc(w_0,\dots w_{k_{n}-1})$, for some $w_i\in \mcw_n$. Since $m_0$ does not occur in the first occurrence of $w_0$ in $w$, this first occurrence is at an interval $[a_0, b_0]$ with $b_0<m_0$. Since some $\sigma\in \Sigma$ must occur in $w_0$ we have the $m_0$ as required. 

From this we see that if $s\in S$ and $\{i:s(i)\notin {b,e}\}$ is bounded below in $\mathbb Z$, then $s$ must have a left tail consisting of the letter $b$, and the similar statement holds for $e$'s and right tails.

 Let $\nu$ be an  shift-invariant measure concentrating on $S$ and $\la \nu_i:i\in I\ra$ be $\nu$'s ergodic decomposition.  Let $T$ be the set of $s\in S$ such that 
there for some $k\in \poZ, s\rest (-\infty, k)$ is constantly $b$ (i.e. those elements of $s$ that are constantly $b$ on a tail going left). Then $T$ is a shift-invariant set. We claim $\nu(T)=0$. If $\nu(T)\ne 0$, then for some $i, \nu_i(T)\ne 0$. Thus without loss of generality we can assume that $\nu$ is ergodic. 
 The ergodic theorem applied to the basic open set 
$\la b\ra$ centered at 0 and the averages ${1\over N}\sum_{-N+1}^0sh^n(\la b\ra)$ shows that 
$\nu(\la b\ra)=1$. The shift invariance and countable completeness of $\nu$ implies that $\nu$ gives the constant $b$ 
sequence measure one, contradicting the assumption that $\nu$ concentrates on $S$.  The proof that $\nu$ 
gives the collection of $s\in S$ with a positive tail constantly $e$ measure zero is very similar, using the 
ergodic averages in the positive direction.

For the reverse implication we show that the collection of $s\in \bk$ such that $\{i:s(i)\notin \{b, e\}\}$ is 
unbounded in both $\poZ^-$ and $\poZ^+$ is a subset of $S$. Let $s$ have this property. Suppose that   we are given  $a, b\in \nn$. We must find an $a'>a$ and a $b'>b$ such that $s\rest[-a',b')\in \mcw_n$.
Choose an $i>a, i'>b$ such that $s(-i)\notin \{b,e\}$ and $s(i')\notin \{b, e\}$. Let $n$ be so large that 
$q_n>i+i'+1$ and consider $s\rest[-10q_{n+1}, 10q_{n+1}]$. This must be a subword of some word $w^*$ in 
$\mcw_{m+1}$ with $m\ge n+1$. Suppose that $w^*=\mcc(w_0, \dots, w_{k_m-1})$.  Since the connected segments of the boundary of $w^*$ are of length $q_{m+1}>i+i'+1$, and neither $-i$ or $i'$ are in a position in $s$ corresponding to the boundary of $w^*$, the positions of $w^*$ corresponding to $-i$ and $i^*$ must be in a segment of $w^*$ of the form $w_k^{l_m-1}$. If $-i-1$ and $i'+1$ are positions in $w^*$ in different copies of $w_k$, then they must be separated by a segment of $b$'s and $e$'s of length $q_m$, a contradiction. Hence they lie in a particular copy of $w_k$. Letting $-a'$ be the beginning position of that copy and $b'$ be the end position,  we have finished the proof of the first claim.

To see the second item, let $\nu$ be an ergodic non-atomic measure. Then, as in the first claim, $\nu$ gives measure zero to those elements of $\bk$ that are constant on a tail in either direction. Hence $\nu$ concentrates on those $s$ that have arbitrarily large positive and negative $i$ in both directions where $s(i)$ is not in $\{b, e\}$. Hence $\nu$ concentrates on $S$. 

If $\nu$ is an arbitrary non-atomic measure then the measures in its ergodic decomposition have to give measure 
zero to those elements of $\bk$ that are constant on a tail in either direction, and hence $\nu$ concentrates  on $S$. 
Thus the last assertion follows from Lemma \ref{unique ergodicity}.\qed

{We now define a canonical factor of a circular system.}
\begin{definition}\label{first appearance of circle factor}
Let $\la k_n, l_n:n\in\nn\ra$ be a coefficient sequence for a circular system with $\sum 1/l_n<\infty$. Let $\Sigma_0=\{*\}$. Define a uniform  circular construction sequence such that each 
$\mcw_n$ has a unique element as follows:
\begin{enumerate}
\item $\mcw_0=\{*\}$ and
\item If $\mcw_n=\{w_n\}$ then $\mcw_{n+1}=\{\mcc(w_n, w_n, \dots w_n)\}$.

\end{enumerate}
Let $\mck$ be the resulting circular system. 
\end{definition}

Let $\bk$ be an arbitrary circular system with  coefficients $\la k_n, l_n\ra$, $\sum 1/l_n<\infty$.  Then 
$\bk$ has a canonical factor isomorphic to $\mck$ as we see by defining the following function:

\begin{equation}\label{definition of factor map}
\pi(x)(i) = \left\{ \begin{array}{ll}x(i) &
\mbox{if $x(i)\in \{b, e\}$} \\
* &\mbox{otherwise}
\end{array}\right.
\end{equation}

The following easy lemma  justifies the terminology of Definition \ref{first appearance of circle factor}:

\begin{lemma}\label{canonical rotation factor} Let $\pi$ be defined by equation \ref{definition of factor map}. Then:
\begin{enumerate}
\item $\pi:\bk\to \mck$ is a Lipshitz  map,
\item $\pi(sh^{\pm 1}(x))=sh^{\pm 1}(\pi(x))$ and thus
\item $\pi$ is a factor map of $\bk$ to $\mck$ and $\bk^{-1}$ to $\mck^{-1}$
\end{enumerate}
\end{lemma}

\begin{definition} We will call $\mck$ the \emph{circle factor}  or \emph{rotation factor} of any circular system with the same construction coefficients $\la k_n, l_n:n\in \nn\ra$. 
\end{definition}
Let $p_n$ and $q_n$ be defined as in equations \ref{pn and qn}, $\alpha_n=p_n/q_n$. Then, if $\sum 1/l_n<\infty$, the sequence of $\alpha_n$ converge to an irrational $\alpha$. Theorem \ref{mck} says that $\mck$ is isomorphic to a rotation by $\alpha$. In the smooth realization of a circular system $\bk$, the factor $\mck$ corresponds a rotation of  the equator.

Because the rotation by $\alpha$ is discrete spectrum, $\mck$ is a factor of the Kronecker factor of a circular system $\bk$. In general it isn't the whole Kronecker factor; in a sequel to this paper we show that if the sequences of words $w_i$ used by $\mcc$ satisfy randomness assumptions, then $\mck$ coincides with the Kronecker factor.

%

\section{Periodic Processes}\label{cyclical approximations}


An important tool in the Anosov-Katok method of conjugacy is the use of periodic processes to build a 
measure preserving transformation. We give a somewhat simplified version  of the method here--it is 
described in more generality in \cite{katoks_book} and \cite{Katok_Steppin}.

 The idea behind the method is the following standard proposition (\cite{Halmos}):
\begin{prop}
Let $T$ be an ergodic measure preserving transformation, $n\in \nn$ and $\epsilon>0$. Then there is a periodic transformation $S$ with period $n$ such that $d_U(S,T)<\epsilon$. In particular $T$ is a weak limit of periodic transformations. 
\end{prop}

Taking this as a starting point we now  describe  a method for determining an ergodic transformation by using a sequence of periodic transformations. Rather than view our transformations as point-maps, we will view a periodic transformation as a periodic permutation of a partition. Our view is less general than that in \cite{Katok_Steppin} in that we take the cycles of the permutation to all be of the same length. Adapting to the general case is routine.

%
%

\begin{definition}
Let $\mcp$ be a partition of the measure space $X$. A \emph{periodic process} is a permutation of the atoms of $\mcp$ such that:
\begin{enumerate}
\item each cycle has the same length,
\item the atoms in each cycle have the same measure.
\end{enumerate}
If all of the atoms of $\mcp$ have the same measure we will call $\mcp$ a \emph{uniform} periodic process.
\end{definition}

It is convenient to view the cycles $\mct_1, \mct_2, \dots \mct_n$ of $\tau$  as ``towers." This is done by arbitrarily choosing an element $B$ of $\mct_i$ and designating it as the base and viewing the $k^{th}$ level of the tower to be $\tau^k(B)$.

A slightly subtle point is that the periodic process is a map that permutes the \emph{partition} and is not defined 
pointwise. We will frequently manipulate periodic processes ${\tau}$ by  considering measure preserving 
transformations $F$ that permute the partition in the same manner as ${\tau}$ does. We will call such an $F$ a 
\emph{pointwise realization} of ${\tau}$.
\medskip

We need to have a notion of convergence to use periodic processes to determine an ergodic transformation. We use a uniform version of $\epsilon$-refinement.
%
%
%
%

\begin{definition}\label{eps approx}
Let ${\tau}$ be a periodic process defined on $\mcp$ and ${\sigma}$ be a periodic process defined on $\mcq$. We 
will say that ${\sigma}$ \emph{$\epsilon$-approximates ${\tau}$} iff
there are disjoint collections of $\mcq$-atoms $\{S_A:A\in \mcp\}$ and  a set $D\subset X$ of measure less than $
\epsilon$ such that for some choice of bases for the towers of $\tau$:
\begin{enumerate}
\item for each $A\in \mcp$ we have $(\bigcup S_A)\setminus D\subseteq A$

\item If $A\in \mcp$ is not on the top level of a $\tau$-tower and  $B\in S_A$ we have  
${\sigma}(B)\setminus D\subseteq {\tau}(A)$

and 
\item for each tower $\mct_k$ of ${\sigma}$ the measures of the intersections of $X\setminus D$ with each level of 
$\mct_k$ are the same.

\end{enumerate}
\end{definition}

This definition is saying that after removing a set of measure less than $\epsilon$, 
the action of ${\sigma}$ is {subordinate} to the action of ${\tau}$, except on the top level of the towers of $\tau$.
We note that in the definition of $\epsilon$-approximation, the second and third clauses imply that we can view 
${\sigma}$ as a periodic process defined on the set $X$ that is built by a ``cutting and stacking" construction from the 
restriction of ${\tau}$ to $X\setminus D$ using subsets of  $D$ as fill sets. {We make this explicit when we compute symbolic representations of limits of periodic systems in section \ref{abstract reps of periodic}.}

\begin{remark}\label{refining towers}
In the version of the Anosov-Katok construction we consider, the situation is somewhat simpler than the general construction given in Definition \ref{eps approx} in that the levels of the towers of $\tau_{n+1}$ are either subsets of or disjoint from the levels of the towers in $\tau_n$. In this case $\bigcup S_A\subseteq A$.

\end{remark}

We will use periodic processes to build an ergodic transformation in the following manner.
\begin{lemma}\label{convergent approximations}
Let $\la \varepsilon_n:n\in\nn\ra$ be a summable sequence of positive numbers and $\la {\tau}_n:n\in\nn\ra$  be a sequence of periodic processes defined on a  sequence of partitions $\la \mcp_n:n\in\nn\ra$.  
Suppose that
\begin{enumerate} 
\item ${\tau}_{n+1}$ $\varepsilon_n$-approximates ${\tau}_n$

and
\item the sequence $\la \mcp_n:n\in\nn\ra$ $\sigma$-generates the measure algebra.
\end{enumerate}
Then there is a unique transformation $T:X\to X$ such that for all $n_0$
\begin{eqnarray}\label{phony convergence}
\lim_{n\to \infty} \mu(\bigcup_{A\in\mcp_{n_0}}({\tau}_nA\Delta TA)) = 0
\end{eqnarray}
\end{lemma} 

\pf For each $n$, let $F_n$ be a pointwise realization of $\tau_n$.  Using the fact that the sequence of partitions 
generates the $\sigma$-algebra, we can apply Lemma \ref{criterion for weak convergence} to see that the sequence 
$\la F_n:n\in\nn\ra$ converge in the weak topology. If $T$ is the limit then clearly equation \ref{phony convergence} 
holds.\qed

We will call a sequence satisfying the hypothesis of Lemma \ref{convergent approximations} a \emph{convergent sequence of periodic processes}. We note that the proof of Lemma \ref{convergent approximations} shows the following.

\begin{lemma}\label{pointwise cyclical} Let $\la \varepsilon_n:n\in\nn\ra$ be a summable sequence of positive numbers.
Suppose that 
$\la \tau_n:n\in \nn\ra$ is a sequence of periodic processes converging to a measure preserving transformation $T$. 
 Let $\la F_n:n\in \nn\ra$ be an arbitrary sequence of measure preserving, invertible transformations such that for each $n$, 
 $\sum_{A\in \mcp_n}\mu(F_nA\Delta \tau_nA)<\varepsilon_n$. 
Then $\la F_n:n\in \nn\ra$ converges weakly to $T$.
\end{lemma}
\pf Same as Lemma \ref{convergent approximations}.\qed

We will use the next lemma to construct isomorphisms between limits of measure preserving transformations.

\begin{lemma}\label{creating isomorphisms}
Fix a summable sequence of positive numbers $\la \varepsilon_n:n\in\nn\ra$. 
Let $(X,\mu)$ and $(Y,\nu)$ be standard measure spaces and $\la T_n:n\in\nn\ra$, $\la S_n:n\in\nn\ra$ be 
measure preserving transformations of $X$ and $Y$ that converge  in the weak topology to $T$ and $S$ respectively. Suppose that $\la \mcp_n:n\in\nn\ra$ is a decreasing sequence of partitions and $\la \phi_n:n\in\nn\ra$ is a sequence of measure preserving transformations such that
\begin{enumerate}
\item $\phi_n:X\to Y$ is an isomorphism between $T_n$ and $S_n$,
\item the sequences $\la \mcp_n:n\in\nn\ra$ and $\la \phi_n(\mcp_n):n\in\nn\ra$ generate the measure algebras of $X$ and $Y$ respectively,
\item $D_\nu(\phi_{n+1}(\mcp_n), \phi_n(\mcp_n))<\varepsilon_n$.

\end{enumerate}
Then the sequence $\la \phi_n:n\in\nn\ra$ converges in the weak topology to an isomorphism between $T$ and $S$.

\end{lemma}
\pf
Conditions 2 and 3 verify the hypothesis of  Lemma \ref{criterion for weak convergence}. Thus the sequence of 
$\phi_n$ converge.  The proof that the limit is an isomorphism is similar.\qed


\section{The Anosov-Katok method of conjugacy}\label{AK method}

We now give a brief exposition of a method developed by Anosov and Katok (\cite{AK}) for realizing abstract measure preserving systems as $C^\infty$-transformations on the unit disk, the annulus or the two-torus ($\mathbb T^2$).\footnote{As noted earlier, we only describe the \emph{untwisted} case of the Anosov-Katok method.} The main result of \cite{AK} is that there is an ergodic measure preserving diffeomorphism of the unit disk in $\poR^2$ that is isomorphic to an irrational rotation of the circle. An important feature of the construction is that the irrational is Liouvillean; the irrationals for which the construction works all  can approximated rapidly by rational numbers.
The Anosov-Katok method has been simplified and  extended by several people (\cite{FSW, FK}). We discuss only two such results.

 Herman (\cite{FayKrik}) proved that if $T$ is a $C^k$-diffeomorphism ($k\ge 2$) that has Diophantine rotation number on the boundary, then  there are  $T$-invariant closed curves arbitrarily close to the boundary, a property violating ergodicity. In particular if $T$ is an  ergodic area--preserving $C^k$-diffeomorphism 
of the unit disk measure with rotation number $\alpha$ on the boundary,  then $\alpha$ is Liouvillean.

A remarkable converse of this theorem is proved in \cite{FSW}, where it is shown that  if $\alpha$ is an arbitrarily Liouvillean irrational then there is a $C^\infty$-diffeomorphism of the torus (or unit disk or annulus) that is isomorphic to rotation by $\alpha$ and has rotation number $\alpha$.

Given these results, it is natural to ask if behavior opposite to that of irrational rotations can be realized as diffeomorphisms. This question was answered in the original paper of Anosov and Katok (\cite{AK})  where it is shown that there are area preserving diffeomorphisms of the disk (or torus or annulus) that are weakly mixing. This result was extended in \cite{FS} to get weakly mixing diffeomorphisms with arbitrary Liouvillean rotation number. 
 
These results are proved by careful examination of the norms of a convergent sequence of diffeomorphisms. We give up this control in our construction which is closer to the original untwisted Anosov-Katok method.

All of the theorems in this section are  variations of the known results appearing in \cite{katoks_book}. What is new in this paper is a symbolic representation of the untwisted Anosov-Katok systems as uniform circular systems. In applications of the results of this paper we are concerned with building circular systems with intricate combinatorial properties. To get diffeomorphisms with these properties we appeal to the results in the next few sections to see that the circular systems can be realized as measure preserving diffeomorphisms.

\subsection{Abstract Anosov-Katok-method}\label{abstract AK}

The Anosov-Katok construction inductively defines a sequence 
$\tau_n$ of periodic processes that converge to the desired transformation $T$. From a high-level point of view, what is useful to us in later applications  is that this construction allows us to 
insert an abitrary finite amount of information into each $\tau_n$ and still make the construction converge to a diffeomorphism provided that a simple equivariance condition is satisfied.  

We elucidate the method on the torus, for convenience. The techniques are easily modified to give symbolic representations of Anosov-Katok diffeomorphisms of the disk or annulus.

We will present the method in two stages. In the first we follow \cite{katoks_book} very closely, repeating the discussion and using the notation there as far as possible, consistent with our later purposes. That construction has three sequences of numbers $\la k_n:n\in \nn\ra$, $\la l_n:n\in\nn\ra$ and $\la s_n:n\in\nn\ra$ as parameters. The assertion will be that if $\la l_n:n\in\nn\ra$ goes  to infinity fast enough (the size of $l_n$ depends on $\la k_m:m\le n\ra$ and $\la s_m:m\le n+1\ra$), $s_{n+1}\le s_n^{k_n}$ and $s_n$ goes to infinity, then the method creates a sequence of periodic approximations that converges to an ergodic transformation $T$.

Given coefficient sequences $\la k_n\ra$ and $\la l_n \ra$ we can build sequences $\la p_n\ra $ and $\la q_n\ra$ as in equation \ref{pn and qn}, starting with $p_0=0$ and $q_0=1$. 
The rationals 
\begin{eqnarray}
\alpha_n={p_n\over q_n}\label{alpha n}
\end{eqnarray}
will approximate a Liouvillean number $\alpha=\lim \alpha_n$. It is clear that this approximation is very fast once we note that 
\begin{eqnarray*}
{\pnpo\over\qnpo}={\pn\over\qn}+{1\over \qnpo}=\alpha_n+{1\over \qnpo}.
\end{eqnarray*}
Thus if the $l_n$ grow fast enough, $\alpha$ is Liouvillean.

The only condition we put on the sequence of $s_n$'s is that
$s_n$ divides $s_{n+1}$ and that they go to infinity.


What we call the \emph{abstract method} defines a transformation on an abstract measure space $(X,\mcb, \mu)$ by periodic processes.  Since approximation by periodic processes produces a weakly convergent sequence we will get no information about the continuity properties of the limit.

There is an auxilliary space {$\mca=\zoo\times [0,1]$} which provides the combinatorial basis for the approximations.  
For the abstract method we will neglect a set of measure zero and view $\mca$ as $\zoo\times \zoo$. 
This is appropriate because
we aren't concerned about continuity properties. In the smooth case we must worry about the boundary; there the action on the two boundary segments $\zoo\times \{0\}$ and $\zoo\times \{1\}$ of the limit transformation will be identical so, \emph{a fortiori}, we are working on the torus.

 We let $S^1$ act on $\mca$ by ``rotation" on the first coordinate. We will denote this additively, viewing the rotation action on $S^1$ as ``addition mod 1" in the $x$-direction of the unit square. Specifically we identify $[0,1)$ with $S^1$ by the map $x\mapsto e^{2\pi i x}$. Then rotation by $\alpha$ corresponds to addition in the exponent.
Given $\alpha\in \mathbb R$ we denote the ``rotation" of the unit interval determined by $\alpha$ as $\mcr_\alpha$  and the horizontal rotation of the rectangle $\mca$ by $\alpha$ as $\rid{\alpha}$.



For positive  $q, s\in\nn$  we let  \hyperlink{IQ}{$\xi^q_s=\mci_q\otimes \mci_s$}; i.e. the partition  of $\mca$ that has atoms of the form 
\[[i/q, (i+1)/q)\times [j/s, (j+1)/s)\]
with $0\le i<q$ and $0\le j<s$. 
For given sequences $\la k_n \ra, \la l_n \ra$ and $\la s_n\ra$, we simplify notation by setting $\xi_n=\mci_{q_{n}}\otimes \mci_{s_{n}}$. If $Z:\mca\to X$ is an invertible measure preserving transformation, then $Z$ defines  partitions 
of $X$ by setting  $\zeta^q_s=Z\xi^q_s$ and $\zeta_n=Z\xi_n$. We will refer to the rectangle $\hoo{i}{q_{n}}\times \hoo{j}{s_{n}}$ as $R^n_{i,j}$ and call $R^n_{i,j}$  the $(i,j)^{th}$ element of $\xi_n$.

If $\alpha=p/q$ with $p, q$ relatively prime, then the atoms of $\xi^q_s$ are permuted by the action $\rid{\alpha}$ and this permutation has $s$ cycles, each of length $q$. 
Conjugating by $Z$ gives a periodic process $\tau$ defined on $X$ with partition $\zeta^q_s$ that has $s$ towers of length $q$. 
When building periodic processes we often want to view $\tau$ as the permutation of the atoms of $\zeta^q_s$ and not as a pointwise realization of a measure preserving map.


In this paper our periodic processes will be strongly uniform and  have  pointwise realizations of the following form:
\begin{eqnarray*}T_n=Z_n\circ \rid{\alpha_{n}}\circ Z_n^{-1}
\end{eqnarray*}
where the sequence of $\alpha_n$ are defined by equations \ref{pn and qn} and \ref{alpha n} and $Z_n$ is a measure isomorphism between $\mca$ and $X$. For notational simplicity we let $Z_{0}=Z$, where $Z:\mca\to X$ is a fixed isomorphism from $(\mca, \mcb, \lambda)$ to $(X, \mcb, \mu)$. For $n\ge 1$,  $Z_n$ will be of the form 
\begin{eqnarray*}
Z_n=Z\circ h_1\circ \dots \circ h_n
\end{eqnarray*}
where each $h_i$ is a measure preserving transformation of $\mca$ that induces a permutation of the atoms of $\mci_{k_{i-1}q_{i-1}}\otimes\mci_{s_{i}}$.  Thus $Z_n:\mca\to X$ is an invertible measure preserving transformation.
Because $\xi_{n}$ refines $\mci_{k_{n-1}q_{n-1}}\otimes\mci_{s_{n}}$, we can view $h_n$ as permuting the atoms of $\xi_n$.

\begin{definition} Since
 $\rid{\alpha_{n}}$  gives a periodic process with partition $\xi_n$,  the map $T_n$ 
induces a periodic process with partition $\zeta_n$, which we take to be $\tau_n$.
When we want to view $\tau_n$ as a collection of towers we take the bases of $\tau_n$ to be  the sets  
$Z_{n}R^{n}_{0,s}$, for $s<s_n$. 
\end{definition}

To start the inductive construction we let $s_0\ge2$ and take $\tau_0$  to be the  periodic process  based on the partition $\xi_0$ induced by  the action on $\xi_0$ given by $\rid{\alpha_0}$ (which is the identity map). Thus $\tau_0$ has $s_0$ towers of height one. 


 What remains is to describe how to pass from $Z_{n}$ to $Z_{n+1}$. The trick (due to Anosov and Katok) is to note that if $h_{n+1}$ commutes with $\rid{\alpha_{n}}$, then: 
\begin{eqnarray*}
T_{n}&=&Z_{n}\rid{\alpha_n}Z_{n}^{-1}\\
&=&Z_{n}h_{n+1}h_{n+1}^{-1}\rid{\alpha_n}Z_{n}^{-1}\\
&=&Z_{n}h_{n+1}\rid{\alpha_n}h_{n+1}^{-1}Z_{n}^{-1}\\
&=&Z_{n+1}\rid{\alpha_{n}}Z_{n+1}^{-1}
\end{eqnarray*}
Consequently, if $\alpha_{n+1}$ is chosen sufficiently close to $\alpha_{n}$ the map $T_{n+1}$ will permute $\xi_{n}$ very similarly to the way that $T_{n}$ does. It follows that the periodic process 
$\tau_{n+1}$ will be very close to the periodic process $\tau_{n}$. We note that determining those $\epsilon$ for which $\tau_{n+1}$ $\epsilon$-approximates $\tau_{n}$ is independent of the choice of pointwise realizations $T_{n}$ and $T_{n+1}$ as it is a property of the partitions.

Summarizing, if $h_{n+1}$ is chosen so that it commutes with $\rid{\alpha_n}$ and the sequence 
$\la \alpha_n:n\in\nn\ra$ converges fast enough, then the sequence of periodic processes we construct 
converges in the weak topology to an invertible measure preserving system. Moreover, we can guarantee that the 
$\alpha_n$ converge arbitrarily fast by choosing the $l_n$-sequence to grow fast.

\begin{remark}\label{the data}
 The foregoing construction of an invertible measure preserving transformation is determined up to isomorphism by
  the following data:
 \begin{enumerate}
 \item The sequences $\la k_n, l_n, s_n:n\in\nn\ra$
 and
 \item the maps $\la h_n:n\in\nn\ra$.

 \end{enumerate}

 \end{remark}

\bigskip

{\bf\noindent{Constructing  $h_{n+1}$:}}
We will build $h_{n+1}$ of a special form. 
We choose  numbers $s_{n+1}$ and $k_n$  with $s_{n+1}$ and $k_n$ being multiples of $s_{n}$ and $s_n^{k_n}\ge s_{n+1}$. We think of 
$k_n$ and $s_{n+1}$ as very large. (In later applications
 $k_n$ will be chosen after $s_{n+1}$ is and will be 
large enough to satisfy some requirements determined by the law of large numbers.) The map $h_{n+1}$ is taken to be a measure preserving transformation permuting $\xi_{n+1}$  that also induces a 
permutation of $\mci_{k_nq_n}\otimes \mci_{s_{n+1}}$; i.e.  it takes atoms of  
$\mci_{k_nq_n}\otimes \mci_{s_{n+1}}$ to atoms of $\mci_{k_nq_n}\otimes \mci_{s_{n+1}}$.

Our transformations will be \emph{untwisted} in the language of \cite{katoks_book}. 
This means that $h_{n+1}$ induces a permutation of the atoms of $\mci_{k_nq_n}\otimes \mci_{s_{n+1}}$ that lie inside $[0, 1/q_{n})\times \zoo$. 
  Since  $\rid{\alpha_n}$ cyclically permutes the towers of $\xi_{n}$ starting with  bases of the form $R^{n}_{0,s_n}$  and the $h_{n+1}$ commute with $\rid{\alpha_n}$, 
  the $h_{n+1}$'s are determined by what they do on $[0, 1/q_{n})\times \zoo$.

Indeed to define $h_{n+1}$ we start with an arbitrary permutation $p$ of $(\mci_{k_nq_n}\otimes \mci_{s_{n+1}})\rest[0, 1/q_{n})\times \zoo$  and let $h_{n+1}^0$ be any pointwise realization of $p$ that gives a permutation of $\xi_{n+1}\rest[0, 1/q_{n})\times \zoo$. We then extend $h^0_{n+1}$ to a pointwise map on $\mca$ that commutes with $\rid{\alpha_n}$ by ``copying over" equivariantly via $\rid{\alpha_n}$.\footnote{For twisted transformations, one does not have complete freedom in defining the injection of this collection of atoms into $\mci_{k_nq_n}\otimes \mci_{s_{n+1}}$, as there may be obstacles to extending it equivariantly. These restrictions are easy to satisfy, but complicate our construction.}

Having defined the $h_{n+1}$  with the numbers $k_n$ and $s_{n+1}$ so that it commutes with $\rid{\alpha_n}$,  we choose $l_n$ large enough that $\alpha_{n+1}$ is very close to 
$\alpha_n$.\footnote{The limit of the sequence only depends on the periodic processes, and hence the permutation of $\xi_n$ determined by $h_n$, rather than $h_n$ as a pointwise map. It will follow that there are functions $l_n^*(x_0, \dots x_{n-1}, y_0, \dots y_{n}, z_0, \dots z_{n+1})$ such that if $l_n\ge l_n^*(l_0, \dots l_{n-1}, k_0\dots k_n, s_0, \dots s_{n+1})$ for all $n$, then the construction converges}.
 Since $h_{n+1}$ permutes the elements of $\mci_{k_nq_n}\otimes \mci_{s_{n+1}}$ it permutes the elements of $\xi_{n+1}$

\begin{lemma}\label{zeta_n generate}
For each $n$ the partition $\zeta_{n+1}$ refines $\zeta_n$ and the collection of partitions $\{\zeta_n:n\in \nn\}$ generates the measure algebra of $X$. 
\end{lemma}
\pf Since $\la \xi_n:n\in\nn\ra$ is a decreasing sequence of partitions that generates the measure algebra of $\mca$, $\la Z\xi_n:n\in\nn\ra$ is a decreasing sequence that generates the measure algebra of $X$.  Each $h_m$ is a permutation of $\xi_n$ for $m\le n$. Consequently, $Z\xi_n=Z_n\xi_n=\zeta_n$. Hence the $\zeta_n$ are decreasing and generate.\qed
%

{\bf\noindent Ergodicity.} We introduce some requirements to guarantee the  ergodicity of our systems. The requirements we state here are  stronger than necessary,\footnote{Cleverer versions of $Z_n$ (e.g.  \cite{FS}) can be constructed that guarantee ergodicity even if every $s_n$ is equal to 2.  {It is not difficult to check that these have symbolic representations very similar to the ones we are describing here. We use Requirements 1-3 to guarantee ergodicity for circular systems with fast growing $l_n$-sequences.}}
 but easy to verify.  We postpone the proof that these requirements imply ergodicity until Theorem \ref{verifying unique ergodicity} in Section \ref{a symb rep for AK}.


Here are our \hypertarget{requirements 1-3}{requirements}:
\begin{description}\label{descs}
  \item[Requirement 1:] The sequence $s_n$ tends to $\infty$.
  \item[Requirement 2:](Strong Uniformity) \label{rqt 2} For each $R^{n}_{0,j}\in \xi_{n}$ and each $s<s_{n+1}$ we have that the cardinality of
  \begin{equation*}
\{t<k_n:h_{n+1}(\hoo{t}{k_nq_n}\times\hoo{s}{s_{n+1}})\subseteq R^n_{0,j} \}\label{equidistribution}
\end{equation*}  is $k_n/s_n$.
\end{description}

\noindent Given $s<s_{n+1}$ we can associate a $k_{n}$-tuple $(j_0, \dots j_{k_n-1})_s$ so that 
\begin{equation}\label{h's}
h_{n+1}(\hoo{t}{k_nq_n}\times\hoo{s}{s_{n+1}})\subseteq R^{n}_{0,j_t}.
\end{equation}
\begin{description}
  \item[Requirement 3:]  We assume that the map $s\mapsto (j_0, \dots j_{k_n-1})_s$ is one-to-one.
\end{description}

\bfni{Discussion:} There are $k_n$ atoms $a$ of $\mci_{k_nq_n}\otimes \mci_{s_{n+1}}$ that lie in the strip 
$[0, 1/q_n)\times [s/s_n, (s+1)/s_n)$. For each such $a$, $h_{n+1}(a)\subset R^n_{0,j}$ for some $j$. If we assign 
this $j$ to $a$, then we get a sequence of $j$'s of length $k_n$. Requirement 2 says that each $j$ occurs $k_n/s_n$ 
times.
It follows that the proportion of the atoms of $\xi_{n+1}$ contained in any strip $[0,1/q_n) \times \hoo{s}{s_{n+1}}$ whose $h_n$-image is a subset of a given $R^n_{0,j}$ is $k_n/s_n$.
 Requirement 3 says that we get different sequences  of $j$'s for different $s$.  
\medskip

Conversely we have {the following lemma which describes the mechanism for inserting arbitrary finite information into each stage of the Anosov-Katok construction.}
\begin{lemma}\label{building hn from words} Let $w_0, \dots, w_{s_{n+1}-1}\subseteq \{0,1,\dots s_n-1\}^{k_n}$ be words such that each $i$ with $0\le i<s_n$ occurs $k_n/s_n$ times in each $w_j$. Then there is an invertible measure preserving 
 $h_{n+1}$ commuting with $\rid{\alpha_n}$ and inducing a  permutation of $\mci_{k_nq_n}\otimes \mci_{s_{n+1}}$  such that if $j_t$ is the $t^{th}$ letter of $w_s$ then 
\[h_{n+1}(\hoo{t}{k_nq_n}\times\hoo{s}{s_{n+1}})\] is a subset of 
$R^{n}_{0, j_t}$.
\end{lemma}
\pf Construct $h_{n+1}^0$ as follows: The atoms of $\mci_{k_nq_n}\otimes \mci_{s_{n+1}}$ partition each atom 
of $\xi_{n}$ into $k_n({s_{n+1}\over s_n})$ pieces. Each index $i$ occurs in each $w_s$ exactly $k_n/s_n$ times and there are $s_{n+1}$ many words $w_s$.  Hence it is possible to construct a bijection between 
$(\mci_{k_nq_n}\otimes\mci_{s_{n+1}})\rest[0, 1/q_n)\times \hoo{i}{s_n}$ and the occurrences of $i$ in all of the words $w_s$. Ranging over $i$ one gets a bijection $b$ between the numbers occurring in all of the words and 
$(\mci_{k_nq_n}\otimes\mci_{s_{n+1}})\rest [0, 1/q_n)\times \zoo$ such that if $j$ is the $t^{th}$ letter of $w_s$ then $b$ associates $\hoo{t}{k_nq_n}\times\hoo{s}{s_{n+1}}$ with an element of $\mci_{k_nq_n}\otimes\mci_{s_{n+1}}\rest R^n_{0,j}$

We can interpret each $w_s$ as assigning numbers to the atoms of $\mci_{k_nq_n}\otimes \mci_{s_{n+1}}$ that are
 subsets of the strip $[0, 1/q_n)\times \hoo{s}{s_{n+1}}$. Hence $b$ can be interpreted as a permutation of 
$(\mci_{k_nq_n}\otimes\mci_{s_{n+1}}) \rest [0,1/q_n)\times \zoo $. We take $h_{n+1}^0$ to be a pointwise realization of $b$ and extend $h^0_{n+1}$ equivariantly. \qed

 \subsection{Approximating Partition Permutations by Diffeomorphisms}
In this section we prove that any permutation of a matrix of rectangles can be well-approximated by a $C^\infty$-measure preserving transformation. We note that similar results were attained independently in \cite{Avila}, \cite{arbieto/matheus} and earlier in \cite{moser}. With the goal of a self-contained exposition,  we present a proof of the theorem here.  

\begin{theorem}\label{smooth sixteen} Let $n$ horizontal lines and $m$ vertical lines divide $\zoc\times \zoc$ into an array of $mn$ equal size rectangles. Let $\sigma$ be a permutation of the rectangles and $\epsilon>0$.  Then there is a $C^\infty$, invertible, measure preserving transformation $\phi$ of $\zoc\times \zoc$ that is the identity on a neighborhood of the boundary of $\zoc\times \zoc$ such that for a set $L$ of Lebesgue measure at least $1-\epsilon$, for all rectangles $R$:
\begin{eqnarray} \mbox{If $x\in L\cap R$, then $\phi(x)\in \sigma(R)$.}\label{obedience}
\end{eqnarray}
\end{theorem}

If $\sigma$ and $\phi$ satisfy the conclusion of Theorem \ref{smooth sixteen}, we will say that $\sigma$ is \emph{$\epsilon$-approximated} by $\phi$. The collection of permutations $\sigma$ that can be $\epsilon$-approximated for all $\epsilon>0$ is closed under composition.

 We first prove a lemma about vertical and horizontal swaps.

\begin{lemma}\label{swap} Consider $[0, 2]\times [0,1]$. Then  for any $\delta>0$, there is a $C^\infty$-measure preserving transformation $\phi_0$ that  is the identity on a neighborhood of the boundary of $[0,2]\times [0,1]$ and for all but $\epsilon$ measure sends $\zoc\times \zoc$ to $[1,2]\times \zoc$ and vice versa.
\end{lemma}

\pf(Lemma \ref{swap}) Let $D\subseteq \poR^2$ be the disk centered at $(0,0)$ that has area $2-\epsilon/2$  and radius $R=R(\epsilon)$. Let $\gamma>0$ be such that the disc of radius $R-\gamma$ has area $2-\epsilon$. It is a standard result that for any positive 
$\gamma$, we can find a $C^{\infty}$ function $f:[0,R]\to [0,\pi]$ 
 such that $f$ is identically equal to $\pi$ on $[0, R-\gamma]$ and is identically equal to zero in an arbitrarily small neighborhood of $R$ in $[0,R]$.

\begin{figure}[h]
\centering
\includegraphics[height=.3\textheight]{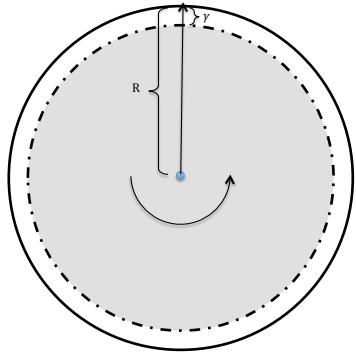}
\caption{The transformation $F$.}
\label{fig1}
\end{figure}

Let $F:D\to D$ be defined in polar coordinates by setting $F(r,\theta)=(r, \theta+f(r))$.
Then $F$ is $C^\infty$, measure preserving, rotates the disk of radius $R-\gamma$ by $\pi$ and is the identity on a neighborhood of the boundary of the disk of radius $R$. 

\begin{figure}[h]
\centering
\includegraphics[height=.25\textheight]{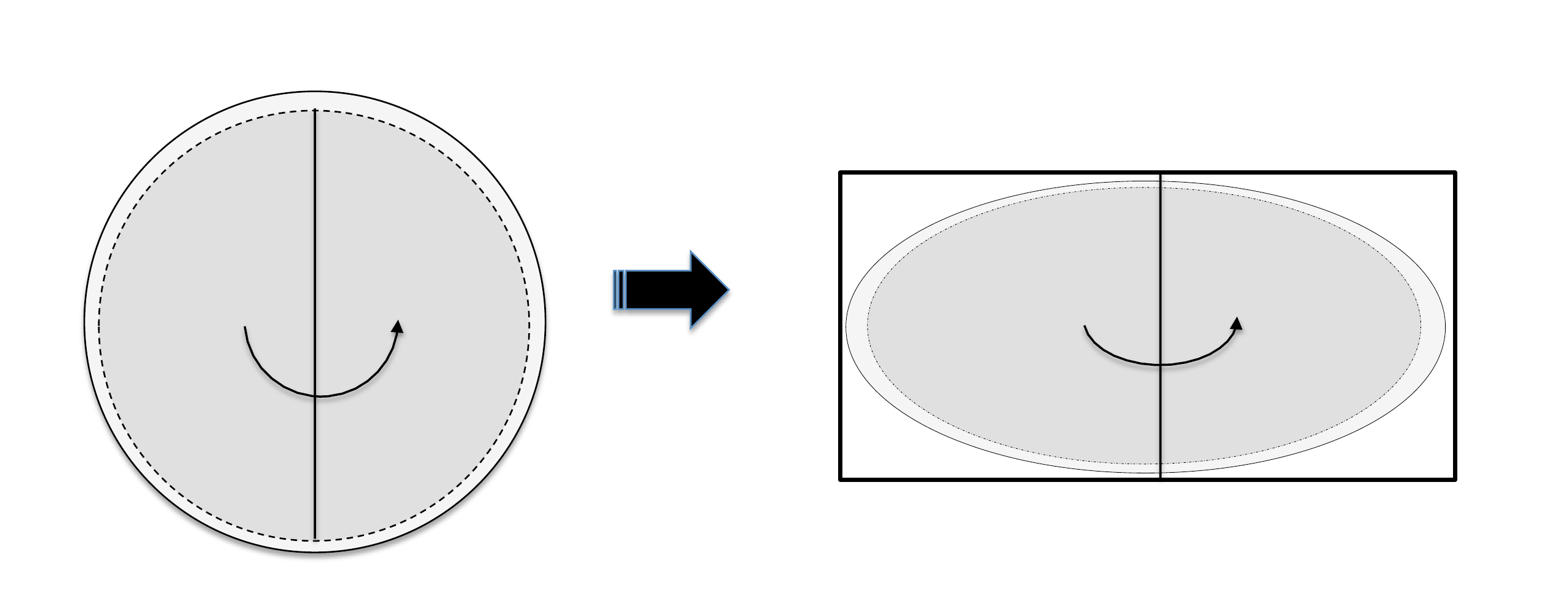}
\caption{The transformation $G$.}
\label{fig2}
\end{figure}

Consider now $[0,2]\times [0,1]$. We can remove a set of measure $\epsilon/2$ near the boundary of $[0,2]\times [0,1]$, so that we are left with a flattened disk $D^*\subseteq [0,2]\times [0,1]$ that has $C^\infty$-boundary. By \cite{moser} there is an invertible measure preserving  
$C^\infty$ function $G:D\to D^*$ that takes the left half disk to $D^*\cap ([0,1]\times [0,1])$ and 
the right half disk to $D^*\cap ([1,2]\times [0,1])$. (See figure \ref{fig2}.)

Then $\phi_0=GFG^{-1}$ is the desired function.
\qed

Clearly Lemma \ref{swap}  can be rescaled arbitrarily. The  result analogous to Lemma \ref{swap} holds for vertically exchanging the two rectangles $\zoc\times [0,1]$ and $\zoc\times [1,2]$  
\bigskip

\pf We now prove Theorem \ref{smooth sixteen}. Number the $mn$ rectanges $\{0, \dots mn-1\}$ in a zig-zag fashion by giving the first row the numbers $0, 1, \dots m-1$ from left to right, enumerating second row from right to left the numbers  $m, m+1, \dots  \dots 2m-1$, the third row from left to right with the numbers $2m, 2m+1 \dots 3m-1$ and so on.

\begin{figure}[h]
\centering
\includegraphics[height=.15\textheight]{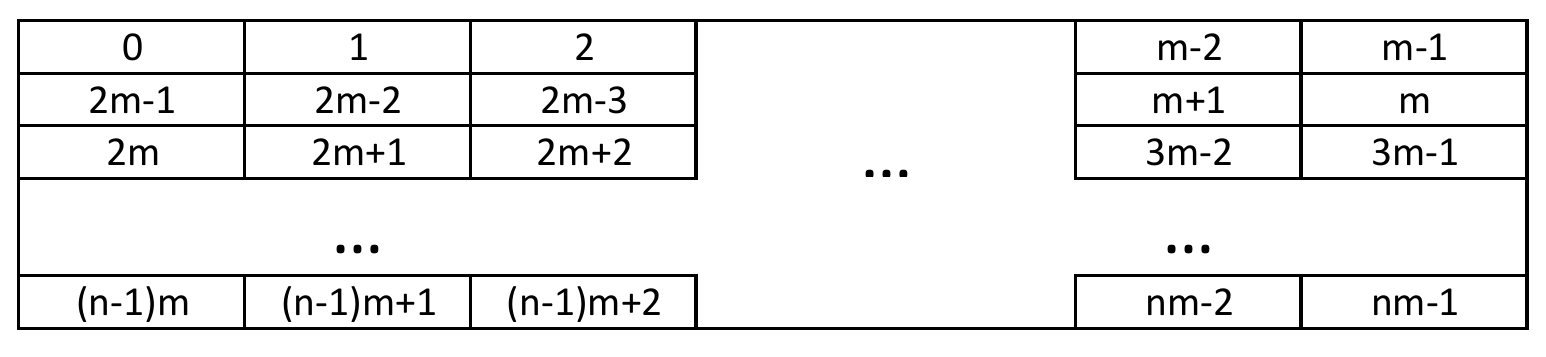}
\caption{Labeling the partition of $[0,1]\times[0,1]$ into rectangles.}
\label{fig3}
\end{figure}


This allows us to view $\sigma$ as a permutation of $\{0, 1, 2, \dots mn-1\}$.  We will call permutations corresponding to exchanging vertically or horizontally adjacent rectangles \emph{swaps}. With this numbering the swaps include all transpositions of the form $(k, k+1)$ for $0\le k< nm-1$. The   smooth measure preserving approximations to swaps given by Lemma \ref{swap} will be called  \emph{$\delta$-approximate swaps}.

Lemma \ref{swap}  implies that for all $k$ with $0\le k<mn-1$ and all $\delta>0$ there is a $C^\infty$ $\delta$-approximate swap $\phi_k$ of the rectangles labelled $k$ and $k+1$. This can be extended to a measure preserving diffeomorphism of $[0,1]\times [0,1]$ by taking $\phi_k$ to be the identity outside the rectangles.

Since every permutation of $mn$ can be written as a composition of less than or equal to $(nm)^2$ transpositions of the form $(k,k+1)$, given any $\sigma$ we can build $\phi$ by taking $\delta$ small enough  and composing $\delta$-approximate swaps corresponding to the transpositions composed to create $\sigma$.\qed

\subsection{Smooth Anosov-Katok-method}\label{smooth AK} 
 We now show how  the Anosov-Katok method of Section \ref{abstract AK} can be used to construct smooth 
 transformations 
 {isomorphic to the abstract transformations constructed in  section \ref{abstract AK}.}
 In this case the measure space $X$ will also be $\mca$, the initial map $Z$ will be the identity,
  and the spatial maps $h_n$ will be replaced by measure preserving diffeomorphisms $h^s_n$ that closely 
  approximate them. Thus for $n\ge 1$, we replace the $Z_n$ by functions 
	  \begin{equation*}
	H_n  =   h^s_1\circ \dots\circ h^s_n.
	\end{equation*}
where the $h^s_i$ are $C^\infty$-measure preserving transformations of $\mca$. We let 
	\[S_n=H_n\rid{\alpha_{n}}H_n^{-1}.\]
Because we want the limit of the $S_n$'s to be a diffeomorphism, we must consider the pointwise properties of the $S_n$; in particular we are no longer working with periodic processes, but concrete realizations of measure preserving transformations. For the limit of the sequence $\la S_n:n\in\nn\ra$ to be a smooth transformation it suffices to arrange that $\|S_{n+1}-S_n\|_{C^{n}}<\varepsilon_n$ for some summable sequence $\varepsilon_n$. Alternately, we fix a metric $d$ giving the $C^\infty$ topology, we can require that $d(S_{n+1}, S_n)<\varepsilon_n$.   These, in turn,  can be arranged by taking $\alpha_{n+1}$ sufficiently close to $\alpha_n$ which is done by choosing  the number $l_n$  large enough.

{\begin{remark}
The transformations $H_n$ will all be equal to the identity on a neighborhood of the boundary of $\mca$. This guarantees that each $S_n$ is equal to $\rid{\alpha_{n+1}}$ on a neighborhood of the boundary of $\mca$. In particular we can view each $S_n$ as a $C^\infty$-measure preserving transformation of the torus. Thus we can view the limit transformation $S$ as a $C^\infty$-measure preserving transformation of the torus that is equal to rotation by $\alpha$ along the line determined by  identifying the top and bottom boundaries of $\mca$. 
\end{remark}}

The main theorem of this section states that we can realize any transformation built by the version of the abstract Anosov-Katok method-of-conjugacy\footnote{
The method of proof Theorem \ref{smoothing things out}
can be easily adapted to other treatments.} as a measure preserving diffeomorphism. This theorem is implicit in the results in \cite{katoks_book}.

\begin{theorem}\label{smoothing things out}

Suppose that $T:\mca\to \mca$ is a measure preserving transformation that is built by the abstract Anosov-Katok-method using a parameter sequence $\la k_n, l_n:n\in\nn\ra$ such that  
the sequence of $l_n$ grow \hyperlink{fe}{fast enough}. Then there is a $C^\infty$-measure preserving transformation $S:\mca\to \mca$ that is measure theoretically isomorphic to $T$. 
\end{theorem}

\pf  Fix a summable sequence $\la \varepsilon_n:n\in\nn\ra$. 
Let $\la T_n:n\in\nn\ra$ be a sequence of transformations built by the abstract Anosov-Katok-method using $\la h_n:n\in \nn\ra$ (thus $T_n=Z_n\rid{\alpha_{n}}Z_n^{-1}$) and let $T$ be the limit of the $T_n$ in the weak topology.  Using Theorem \ref{smooth sixteen} we will build smooth approximations 
$\la h_n^s:n\in\nn\ra$ and define a sequence of smooth transformations by setting $H_0$ be the identity map and for $n\ge 1$, setting 
$H_n=h^s_1\circ \dots h^s_n$. We show that if the $h^s_n$ approximate $h_n$ closely enough and $\la\alpha_n:n\in\nn\ra$ converges fast enough, then the sequence $S_n=H_n\rid{\alpha_{n}}H_n^{-1}$ converges to a $C^\infty$-measure preserving transformation $S$ and that $S$ is isomorphic to  $T$ as a measure preserving transformation of $(\mca, \mcb, \lambda)$.

To see that rapid convergence of  $\la \alpha_n:n\in\nn\ra$  implies that the $S_n$ converge in the $C^\infty$-topology,    we note that  since $h^s_n$ commutes with $\rid{\alpha_n}$  for each $n$ it follows that:
\begin{eqnarray}
S_{n}&=&H_{n}\rid{\alpha_n}H_{n}^{-1}\notag\\
	&=&H_{n}h_{n+1}^s\rid{\alpha_n}(h_{n+1}^s)^{-1}H_{n}^{-1}.\notag
\end{eqnarray}
Hence, by the continuity of $d^\infty$ with respect to composition, if we take $\alpha_{n+1}$ close enough to $\alpha_n$, we can arrange that $d^\infty(S_{n+1},S_{n})<\varepsilon_n$.

Let $n\ge 0$.
 To define $h_{n+1}^s$  we need to both approximate $h_{n+1}$ and make $h_{n+1}^s$ commute with $\rid{\alpha_n}$. We use Theorem \ref{smooth sixteen} to choose a smooth measure preserving, invertible $h'_{n+1}:[0,1/q_{n}]\times [0,1]\to[0,1/q_{n}]\times [0,1]$ such that: 
\begin{enumerate}
\item If $\sigma$ is the permutation of the atoms of  $\mci_{k_{n}q_{n}}\otimes \mci_{s_{n}}$ inside $[0,1/q_{n})\times [0,1)$ determined by $h_{n+1}$, then  for all but a set of measure $\epsilon_{n+1}/q_{n}$, we satisfy equation \ref{obedience}; i.e. for the atoms $R$ of $\mci_{k_{n}q_{n}}\otimes \mci_{s_{n}}$ inside $[0,1/q_{n})\times [0,1)$, we know that the vast majority of $x\in R$ have $h_{n+1}'(x)\in \sigma(R)$,

and 
\item  the function $h'_{n+1}$ is the identity map on a neighborhood of the boundary of $[0,1/q_{n}]\times [0,1]$
\end{enumerate}

Since $h_{n+1}'$ is the identity on a neighborhood of the boundary of $[0,1/q_{n}]\times [0,1]$, we can copy it onto each $[i/q_n, (i+1)/q_n]\times [0,1]$ and thereby extend it to a $C^\infty$-measure preserving $h^s_{n+1}$ that commutes with $\rid{\alpha_{n}}$. If $\sigma$ is the $\rid{\alpha_{n}}$-equivariant permutation of $\xi_{n+1}$ determined by $h_{n+1}$,  then there is a set $L_{n+1}$ of measure  at least $1-\epsilon_{n+1}$ such that for all atoms $R$ of $\xi_{n+1}$ and $x\in L_{n+1}$:
\begin{eqnarray}\label{Ln}
\mbox{If $x\in L_{n+1}\cap R$, then $h_{n+1}^s(x)\in \sigma(R)$.}
\end{eqnarray}

We now use Lemma \ref{creating isomorphisms} to show that $T$ is isomorphic to $S$.
We let $\phi_n=H_n\circ Z_n^{-1}$ and $\mcp_n=\zeta_n$. 
\[\begin{diagram}
\node{}	 \node{\mathcal A}	\arrow{sw,r}{Z_n}\arrow{se,r}{H_n}  \node{}\\
\node{X}		\arrow[2]{e,r}{\phi_n}\node{}	\node{\mca}
\end{diagram}
\]
We must verify conditions 1)-3) of Lemma \ref{creating isomorphisms}. It is clear that $\phi_n$ is an isomorphism between $T_n$ and $S_n$ since $Z_n$ and $H_n$ are isomorphisms between $\rid{\alpha_{n+1}}$ and $T_n$ and $S_n$ respectively.  Thus condition 1) is clear. To see condition 2) we must prove the following claim:



\medskip

\noindent{\bf Claim} The partition $\la\zeta_n:n\in\nn\ra$ generates the measure algebra of $X$ and   $\la\phi_n(\zeta_n):n\in\nn\ra$ generates the measure algebra of $\mca$.
\smallskip

\pf (Claim) That the $\zeta_n$'s generate is the content 
of Lemma \ref{zeta_n generate}. We must check that the $\phi_n(\zeta_n)$ generate. Since $\phi_n=Z_n^{-1}H_n$, this is equivalent to the statement that the $H_n\xi_n$ generate.


Since the sequence $\varepsilon_n$ is summable, the Borel-Cantelli lemma tells us that there is an increasing sequence of sets $\la G_n:n\in\nn\ra$ with $G_n\subseteq \mca$  such that the measures of  ${\lambda}(G_n)$ approach $1$ and for all $m\ge n$ $h_m$ and $h_m^s$ permute that partition $\xi_m\rest G_n$ the same way. 

Explicitly: for all $\epsilon>0$ we can find an $n$ so large that:
\begin{eqnarray}
G_n=_{def}L_n\cap\bigcap_{m> n}(h^s_{n+1}\circ h^s_{n+2}\circ \dots h^s_{m})^{-1}L_m\label{Gn}
\end{eqnarray}
has measure at least $1-\epsilon$.\footnote{The $L_n$'s are defined in formula \ref{Ln}.} By definition, for all $x\in G_n$ and all $m>n, h_{n+1}\circ\dots \circ h_{m}(x)\in L_m$.

Fix a measurable set $D\subseteq \mca$ and a $\delta>0$. We must find a large enough $n$ that the atoms of $H_n\xi_n$ can be used to approximate $D$ within a set of measure  $\delta$. We first choose an $n_0$ so large that ${\lambda}(G_{n_0})>1-\delta/100$. Let $D'=H_{n_0}^{-1}(D)$.

Since the $\xi_n$'s generate, we can find an $n\ge n_0$ and a collection $C'$ of atoms of $\xi_{n}$ such that $(\bigcup C') \Delta D'$ has measure less than $\delta/100$.  Hence we can find a collection of atoms of $H_{n_0}\xi_{n}$ whose union approximates $D$ within $\delta/100$.

We note that for $n_0<m\le n$, $h_m$ permutes the atoms of $\xi_n$. From the definition of $G_{n_0}$, for $n_0<m\le n, a\in \xi_n$ and 
$x\in a\cap G_{n_0}$, we know that  $h_{n_0+1}^s\circ h_{n_0+2}^s\circ \dots h^s_m(x)$ belongs to the same atom of $\xi$ as $h_{n_0+1}\circ \dots h_m(x)$ does. It follows that 
\begin{eqnarray} {\lambda}(\bigcup_{a\in \xi_n}(h_{n_0+1}^s\circ h_{n_0+2}^s\circ \dots h^s_m(a)\ \Delta\  h_{n_0+1}\circ \dots h_m(a)))<\delta/50\notag
\end{eqnarray}


Since $h_{n_0+1}, \dots h_n$ permute the atoms of $\xi_n$, we can define a bijection  $\sigma:\xi_n\to \xi_n$ such that 

\begin{eqnarray*} {\lambda}(\bigcup_{a\in \xi_n}(h_{n_0+1}^s\circ h_{n_0+2}^s\circ \dots h^s_n(a)\ \Delta\  \sigma(a))<\delta/50
\end{eqnarray*}
Thus:
\begin{eqnarray*}{\lambda}(\bigcup_{a\in \xi_n}[H_n(a)\ \Delta\  H_{n_0}(\sigma(a))])<\delta/50,
\end{eqnarray*}
or equivalently

\begin{eqnarray*}{\lambda}(\bigcup_{a\in \xi_n}[H_n(\sigma^{-1}(a))\ \Delta\  H_{n_0}(a)])<\delta/50,
\end{eqnarray*}

Since $D$ can be approximated within $\delta/100$ by a union of atoms of $H_{n_0}\xi_n$ it can be approximated within $\delta/25$ by a union of atoms of $H_n(\xi_n)$. We have verified the claim.\qed

 To use Lemma \ref {creating isomorphisms}, we are left with showing  that $\phi_{n+1}(\zeta_n)$ $\varepsilon_n$-approximates $\phi_n(\zeta_n)$.
Chasing the following diagram, where the leftmost and rightmost triangles commute:
\[\begin{diagram}
\node{}	 	\node{} \node{\mca}\arrow{wsw,l}{Z_{n+1}}\arrow{sw,r}{h_{n+1}}\arrow{se,r}{h_{n+1}^s} \arrow{ese,l}{H_{n+1}}	\node{}  \node{}\\
\node{X}		\node{\mca}\arrow{w,r}{Z_n}\arrow[2]{e,r}{id} \node{} \node{\mca}\arrow[2]{w,e}{}\arrow{e,r}{H_n}	\node{\mca}
\end{diagram}
\]
we see that
\begin{eqnarray*}
\phi_n(\zeta_n)	&=	&H_nZ_n^{-1}(\zeta_n)\\
	&	=&H_nh_{n+1}Z_{n+1}^{-1}(\zeta_n) 
\end{eqnarray*}
while
\begin{eqnarray*}
\phi_{n+1}(\zeta_n)	&	=&H_{n+1}Z_{n+1}^{-1}(\zeta_n)\\
	&	=&H_nh^s_{n+1}Z_{n+1}^{-1}(\zeta_n).
\end{eqnarray*}
Letting $\mcq_n=Z_{n+1}^{-1}(\zeta_n)$, we need to see that $H_nh^s_{n+1}(\mcq_n)$ $\varepsilon_{n}$-approximates $H_nh_{n+1}(\mcq_n)$. This follows easily, since $H_n$ is measure preserving and $h^s_{n+1}$ was chosen so that $h^s_{n+1}(\mcq_n)$ closely approximated $h_{n+1}(\mcq_n)$.\qed
\begin{remark}\label{going to a metric}
Let $d^\infty$ be a Polish metric inducing the $C^\infty$-topology, then by taking
 $\alpha_{n+1}$ is sufficiently close to $\alpha_n$, we can arrange $d^\infty(S_{n+1},S_n)<\varepsilon_n/4$. 
\end{remark}
\hypertarget{fe}{{\noindent{\bf What does \emph{fast enough} mean?}}} 
  We do not produce explicit lower bounds on the speed of growth of the $l_n$'s; instead we give inductive lower bounds. Here is what we mean by \emph{fast enough} in the statement  of Theorem \ref{smoothing things out}.  

Fix a summable sequence $\la \varepsilon_n:n\in\nn\ra$; without loss of generality 
\[\varepsilon_n/4>\sum_{m>n}\varepsilon_m.\]  Fix a metric $d^\infty$ inducing the $C^\infty$-topology on the diffeomorphisms of the 2-torus.
For each choice $\la k_i:i\le n\ra, \la l_i:i<n\ra$, and  $\la s_i:i\le n+1\ra$ there are only finitely many permutations of the relevant partitions and thus only finitely many choices for the periodic processes determined by $\la h_i:i\le n\ra$.  Hence 
there is a single number $l^*_n=l^*_n(\la k_i:i\le n\ra, \la l_i:i<n\ra,\la s_i:i\le n+1\ra)$ such that for all $l_n\ge  l_n^*$ we can choose a smooth approximation $h_n^s$ such that:
 \begin{equation}\label{key to cont}
 d^\infty(S_n,S_{n+1})<\varepsilon_n/4.
 \end{equation} 
 Without loss of generality we can assume that for all $n, l^*_n>2^{n}$.

 Indeed we can say more. In our construction we have the inequality that $s_{n+1}\le s_n^{k_n}$. Hence there is a sequence of  bounds 
$b_n=b_n(|\Sigma|, \la k_{i}:i<n\ra)$ such that $s_n\le b_n$. By ranging over all  $s_n\le b_n$ in the previous paragraph we get a sequence $l_n^*=l_n^*(\la k_i:i<n+1\ra,\la l_i:i<n\ra)$ such that any sequence $l_n\ge l_n^*$ grows fast enough for the hypothesis of Theorem \ref{smoothing things out}.

From now on by  \emph{fast enough} we mean a sequence of $l_n$ with $l_n\ge l^*_n$.

\section{A symbolic representation of Anosov-Katok-systems}\label{symbolic representations of AK}

In this section we show two theorems. The first (Theorem \ref{circular as smooth}), which is used in the sequels, is that if $\bk$ is a strongly uniform circular system with fast growing coefficients then $\bk$ is isomorphic to a measure preserving Anosov-Katok diffeomorphism of $\mca$. As usual $\mca$ is our proxy for the unit disk, annulus or $\bt^2$.

The second theorem (Theorem \ref{verifying unique ergodicity}), which we hope is of independent interest, is that  if $T$ is built by the (untwisted) Anosov-Katok method with coefficients $\la k_n, l_n, s_n:n\in\nn\ra$ and $\la l_n:n\in\nn\ra$ grows fast enough then $T$ has a  representation as a circular symbolic system with the same coefficients. 

By Theorem \ref{smoothing things out}, to represent the Anosov-Katok diffeomorphisms it suffices to represent the transformations built by the abstract method. We will use the notation (such as $k_n, l_n, Z_n, X$) from that section (\ref{abstract AK}). We take $X=\mca$ and $Z_0=Z$ to be the identity map.

\subsection{Symbolic representations of  periodic processes.}\label{abstract reps of periodic} 

We begin with a very general discussion of how periodic processes can be viewed as symbolic systems. Our symbolic description is a variant of standard cutting and stacking constructions where spacers are added at both the bottom (the ``beginning," denoted by $b$'s) and at the top (the ``end," denoted by $e$'s). Our representation is given explicitly for the case where all of the cycles in the periodic process have the same length. It is straightforward to adapt our analysis to the general case.

  A sequence of periodic processes converging sufficiently rapidly to a transformation $T$ gives rise to a symbolic presentation of  a factor of $T$ that has  a special form. Here is how this works. Let $\la \varepsilon_n:n\in\nn\ra$ be a summable sequence.
Fix a sequence $\la {\tau}_n:n\in\nn\ra$ of periodic processes converging to a transformation $T$ with partitions $\la \mcp_n:n\in\nn\ra$. Suppose that the length of the towers corresponding to ${\tau}_n$ is $q_n$ and that $\tau_{n+1}$ $\varepsilon_n$-approximates $\tau_n$.  
Since $\tau_{n+1}$ $\epsilon_n$-approximates $\tau_n$, we can find a set $D_n$ (as in Definition \ref{eps approx}) such that on the complement of $D_n$, the action of $\tau_{n+1}$ is subordinate to $\tau_n$.  By skipping a finite number of steps at the beginning of the approximation we can assume that $\mu(\bigcup D_n)<1/2$ and 
$q_0=1$. 

Setting $G_n=X\setminus \bigcup_{m\ge n} D_n$ we get an increasing sequence of sets $\la G_n:n\in\nn\ra$ such that:
\begin{enumerate} 
\item the measure  of $G_n$ goes to one, 
\item on $G_n$, the sequence $\la \mcp_m:m\ge n\ra$ is a decreasing sequence of partitions

and 
\item for each tower $\mct_k$ of ${\tau}_n$, the measures of the intersections of $G_n$ with the levels of 
$\mct_k$ are all the same.
\end{enumerate}

We describe a  sequence of collections of sets $\mcq_0$, $\la B_n:n\in \nn\ra$ and $\la E_n:n\in\nn\ra$ such that 
$\mcq_0\cup\{B_n, E_n\}$ is a partition of $G_n$.\footnote{The use of $B$'s and $E$'s 
correspond correspond to the mnemonic \emph{beginning} and \emph{end}. }

   Let $\la\mct_i:i< s_0\ra$ be the collection of towers for ${\tau}_0$. Because we assume that $q_0=1$, each of 
   the 
   towers of $\tau_0$ has height 1. Let $\mcq_0$ be the partition of $G_0$ into $s_0$ pieces  consisting of the sets in 
   these towers of ${\tau}_0$ intersected with $G_0$. Let $B_0=E_0=\emptyset$.

\begin{figure}[!h]
\centering
\includegraphics[height=.8\textheight]{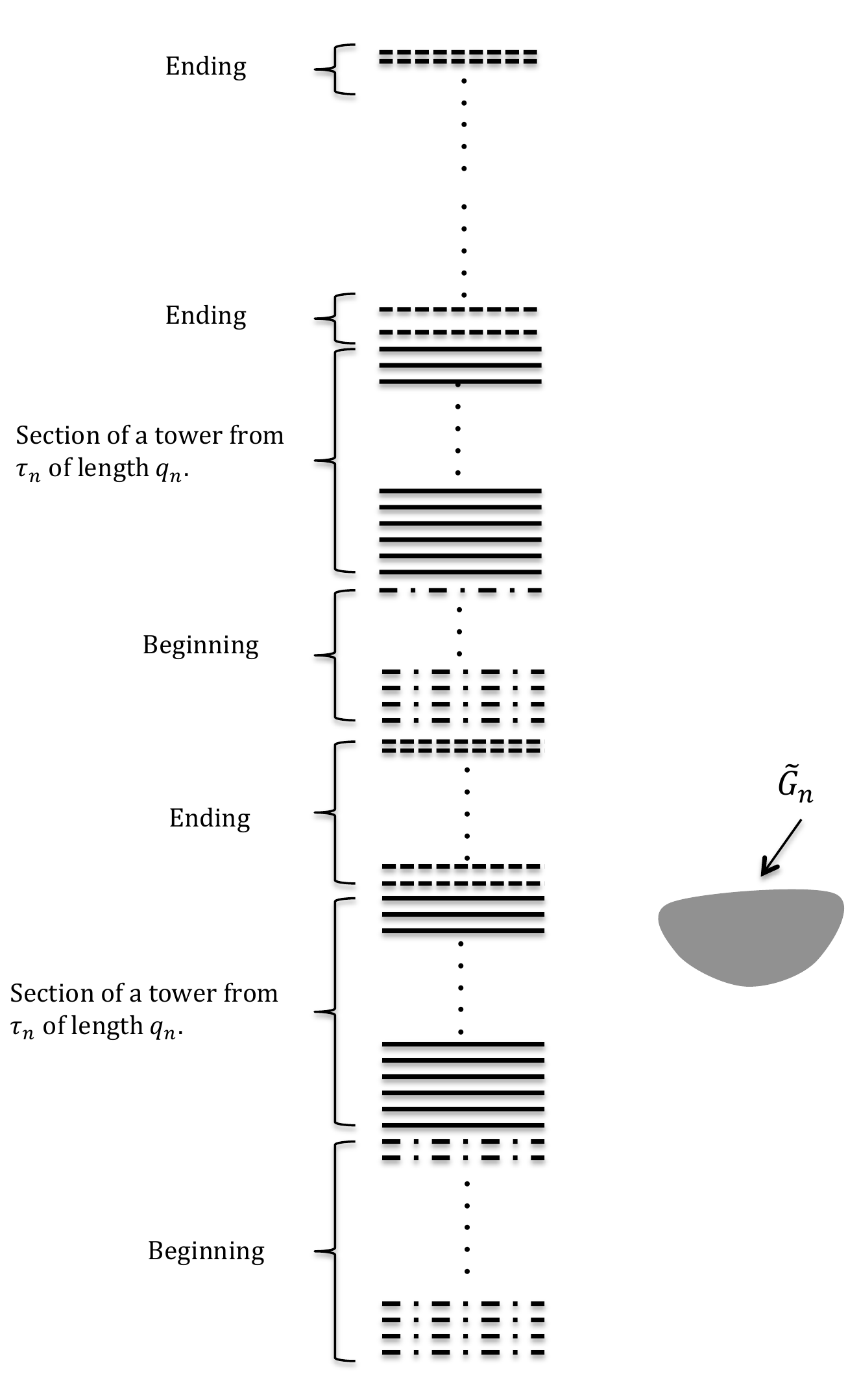}
\caption{A typical tower in $\tau_{n+1}$ built from sections of towers in $\tau_n$ filled in with beginning and ending levels.}
\label{fig4}
\end{figure}

We inductively define the $B_n$ and $E_n$ so that $B_{n+1}\supseteq B_n$ and $E_{n+1}\supseteq E_n$. Suppose that $B_n, E_n$ are defined, and we want to define $B_{n+1}, E_{n+1}$. By assumption each tower of ${\tau}_{n+1}$ restricted to $G_{n+1}$ consists of
\begin{enumerate}
\item contiguous sequences of length $q_n$ consisting of portions of consecutive levels of   the towers of ${\tau}_n$    

interspersed with 
\item new levels of the towers of ${\tau}_{n+1}$ intersected with $G_{n+1}$.
\end{enumerate}
The interspersed levels are first divided into maximal contiguous portions.
 We now arbitrarily divide the levels of each of these portions into two contiguous subcollections of levels.\footnote{If a tower $\mct$ of ${\tau}_{n+1}$ both begins and ends with levels not in a tower of ${\tau}_n$ then we include the new top portion with the new bottom portion  and view it as one contiguous block when we do this division. This is consistent with our view of the ${\tau}_n$ as periodic.}
 One of these subcollections comes before the other in the natural ordering of the tower by the cyclic permutation ${\tau}_{n+1}$. We view the subcollection that comes first as \emph{ending} a block of levels coming from ${\tau}_n$ and the second as \emph{beginning} the next block of levels. The set $B_{n+1}$ consists of the union of $B_n$ with the unions of the points in these second contiguous subcollections of levels and $E_{n+1}$ as $E_n$ together with the first contiguous subcollection of levels (see figure \ref{fig4}).

The partition $\mcq=_{def}\mcq_0\cup \{\bigcup B_n, \bigcup E_n\}$ will generate the transformations we eventually construct, but may not do so in general. We now describe a symbolic representation of the factor of $T$ determined by $\mcq$.

 Suppose that the number of towers in the partition corresponding to ${\tau}_0$ is $s_0$. Let $\Sigma$ be an alphabet of cardinality $s_0$. We view 
 $\Sigma$ as indexing the partition $\mcq_0$; e.g. if $\mcq_0=\{A_i:i\in I\}$ then $\Sigma=\{a_i:i\in I\}$. The alphabet of our symbolic shift will be $\{b, e\}\cup \Sigma$, where $b$ and $e$ are symbols not in $\Sigma$. Proceeding as in Section \ref{presentations}, we define a factor map $\pi:X\to (\Sigma\cup \{b, e\})^\poZ$ by letting:
 \begin{enumerate}
  \item  $\phi(x)(m)=a_i$ iff $T^m(x)\in A_i$,
  \item $\phi(x)(m)=b$ iff $T^m(x)\in \bigcup_n B_n$, and
  \item $\phi(x)(m)=e$ iff $T^m(x)\in \bigcup_n E_n$.
\end{enumerate}
 The symbolic representation is the resulting system {$((\Sigma\cup\{b,e\})^\poZ, \mcc, \phi^*\mu, sh)$}.
 
 We now examine the construction of the approximations and the partition $\mcq_0\cup\{\bigcup B_n, \bigcup E_n\}$ to get a clear description  of the support of the measure $\phi^*\mu$.  We inductively define a sequence of collections of words $\la \mcw_n:n\in \nn\ra$ with the following properties:
 \begin{enumerate}
 \item  There is a surjection  $\phi_n$ from  $\{\mct\cap G_n:\mct$ is a tower for ${\tau}_n$ and $G_n\cap \mct$ is non-empty$\}$ to $\mcw_n$.
 \item The length of each word in $\mcw_n$ is $q_n$.
  \end{enumerate}

 Let $\mcw_0=\Sigma$.
  Suppose that we have defined $\mcw_{n}$  and $\phi_n$
We want to define $\mcw_{n+1}$ and $\phi_{n+1}$. The levels of each $\mcs\cap G_n$ where $\mcs$ is a tower of ${\tau}_{n+1}$  come in blocks of three kinds:
\begin{enumerate}
\item a contiguous block of length $q_n$ coming from a tower of ${\tau}_n$,
\item a contiguous block of levels in $B_{n+1}\setminus B_n$, 

and 
\item a contiguous block of levels in $E_{n+1}\setminus E_n$.
\end{enumerate}

The word in $\mcw_{n+1}$ we associate with $\mcs\cap G_n$ via $\phi_{n+1}$ is the word of length $q_{n+1}$ whose $j^{th}$ letter is $v$ if 
\begin{enumerate}
\item the $j^{th}$ level is in the $k^{th}$ place of a block of  coming from a tower $\mct$ of ${\tau}_n$
and the $k^{th}$ letter of $\phi_n(\mct)$ is $v$,

or

\item $v=b$ and the $j^{th}$ level is in $B_{n+1}\setminus B_n$, 

or
\item $v=e$ and the $j^{th}$ level is in $E_{n+1}\setminus E_n$.

\end{enumerate}

Working as in Section \ref{symbolic shifts}, we let $\bk$ be the collection of $x\in \Sigma^\poZ$ such that every finite contiguous subword of $x$ occurs inside some $w\in \mcw_n$. Then $\bk$ is a closed shift-invariant set that constitutes the support of $\phi^*\mu$.

We will use this technique to represent the smooth transformations we construct. It will follow from the ``Requirements 1)-3)" that in our representation we can choose the partition $\mcq$ so that it generates the transformation $T$ and the system $\bk$ will satisfy Lemma \ref{unique ergodicity}. 
 In particular, the unique non-atomic  measure will be $\phi^*\mu$ and we will have a symbolic representation of our transformation $T$.

\subsection{The dynamical and geometric orderings}
Fix a rational $\alpha={p/q}$ in reduced form. We endow the partition $\mci_q$ of $\zoo$ with two orderings. The first is straightforward: ordering these intervals from left to right according to their left endpoint gives us the \emph{geometric} ordering; in other words, $I^q_i<I^q_j$ iff $i<j$.

  The rotation by $\alpha$ gives us a different ordering which we call the \emph{dynamical} ordering, $<_d$, which we now define explicitly:

We set an interval $ I^q_i<_dI^q_j$ iff there are $n<m<q$ such that $np\equiv i\mod q$ and $mp \equiv j\mod q$.

This can be rephrased conceptually as follows. The first interval in the ordering is $I^q_0$. Repeatedly applying $\mcr_\alpha$,   we get a sequence of intervals $I_0, \mcr_\alpha I_0, \mcr_\alpha^2I_0 \dots,$ $\mcr_\alpha^{q-1}I_0$. The list 
\[\la I_0, \mcr_\alpha I_0, \mcr_\alpha^2I_0 \dots \mcr_\alpha^{q-1}I_0\ra\] gives the dynamical ordering of $\mci_q$.
\begin{remark}\label{the jis}
If we let $j_i\equiv(p)^{-1}i\mod q$ (with $0\le j_i<q$), then the $i^{th}$ interval in the geometric ordering,
 $I^q_i$, is $j_i^{th}$ interval in the dynamical ordering.
\end{remark}

\subsection{Transects}

\begin{quotation}
{\noindent Wikipedia Entry: \emph{A transect is a path along which one records and counts occurrences of the phenomena of study.}\footnote{As of August 15, 2011.} 
}
\end{quotation}

\noindent We now see how the periodic processes $\tau_n$ and $\tau_{n+1}$ compare.

\subsubsection{Without the partitions}\label{wo parts}
Let $\alpha=p/q$ and $\beta=p'/q'$, where we assume $\beta=\alpha +{1\over klq^2}$.\footnote{{Thus $q'=klq^2$.}} We compare how the 
dynamical ordering determined by $\beta$ interacts with the dynamical ordering determined by $\alpha$.  In the discussion that follows we use $j_i$ to denote $(p)^{-1}i \mod q$, i.e. we refer to the dynamical ordering of $\mci_q$ with respect to $\alpha$.

If $J=\hoo{t'}{q'}$ is a subinterval of $I=\hoo{t}{q}$, then $\mcr_\beta J$ is a subinterval of $\mcr_\alpha I$ 
\emph{unless} $J$ is the geometrically last interval in the subdivision of $I$ into intervals of length $q'$. In the latter 
case $\mcr_\beta J$ is the first subinterval of the geometric successor of $\mcr_\alpha I$. 


\begin{figure}[h]
\centering
\includegraphics[height=.20\textheight]{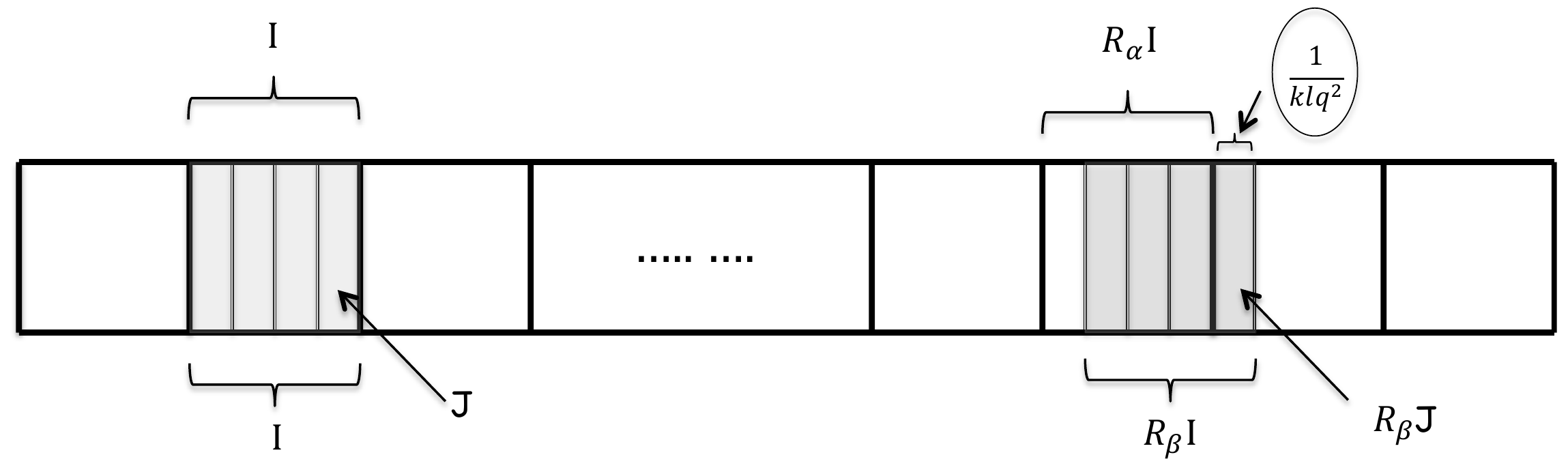}
\caption{A diagram of $\mcr_\alpha$ and $\mcr_\beta$ acting on the first coordinate of $[0,1]\times [0,1]$, showing where $\mcr_\alpha$, $\mcr_\beta$ send $I\times [0,1]$ and $J\times [0,1]$. }
\label{fig5}
\end{figure}

Restating this, if $\mcr_\alpha I$ is the $j_i^{th}$ interval in the dynamical ordering, and $J$ is the geometrically last subinterval of $I$, then $\mcr_\beta J$ is the geometrically first subinterval of the $j_{i+1}^{st}$ interval in the dynamical ordering of $\mci_q$. 

If $J=[0, 1/q')$ is the geometrically first subinterval of  $I=I^q_0$, then 
 for $0\le n<klq$, $\mcr^n_\beta J\subseteq \mcr_\alpha^nI^q_0$ and $\mcr_\beta^{klq}J$ is the geometrically first subinterval of $I^q_1$. Since $\mcr_\alpha^{mq}I=I$ for all $m$, we have $\mcr_\beta^{mq}J$ is a subinterval of $I$ for $0\le m<kl$.
 
%

In general, the $n$'s between $mklq$ and $(m+1)klq$ can be split into three pieces:
\begin{description}
\item[-] the {\emph{beginning}} interval $[mklq, mklq+q-j_m)$ that has length $q-j_m$
\item[-] the \emph{middle} interval $[mklq+q-j_m, mklq+q-j_m+(kl-1)q)$ that has length $(kl-1)q$

and
\item[-] the \emph{end} interval $[mklq+q-j_m+(kl-1)q, (m+1)klq)$ that has length $j_m$.

\end{description}
Assume inductively that $\mcr^{mklq}_\beta J$ is the geometrically first subinterval of $I^q_m$. Then for 
$mklq\le n<mklq+q-j_m$, $\mcr^n_\beta J$ transits the intervals  in places $j_m$, $j_{m}+1$ \dots ${q-1}$ in the 
dynamical  ordering  
of $\mci_q$. The next application of $\mcr_\beta$ puts the orbit in the first interval of $\mci_q$ in the dynamical ordering, which coincides with the first interval in the geometric ordering.

For the \emph{middle} range $n$'s (those $n$ for which $mkl+q-j_m\le n<mkl+q-j_m+(kl-1)q$): {Starting as a subset of 
$I^q_0$,}  $\mcr^n_\beta J$ transits the intervals  $\mci_q$ following the dynamical ordering coming from $\alpha$. 

\begin{figure}[!h]
\centering
\includegraphics[height=.50\textheight]{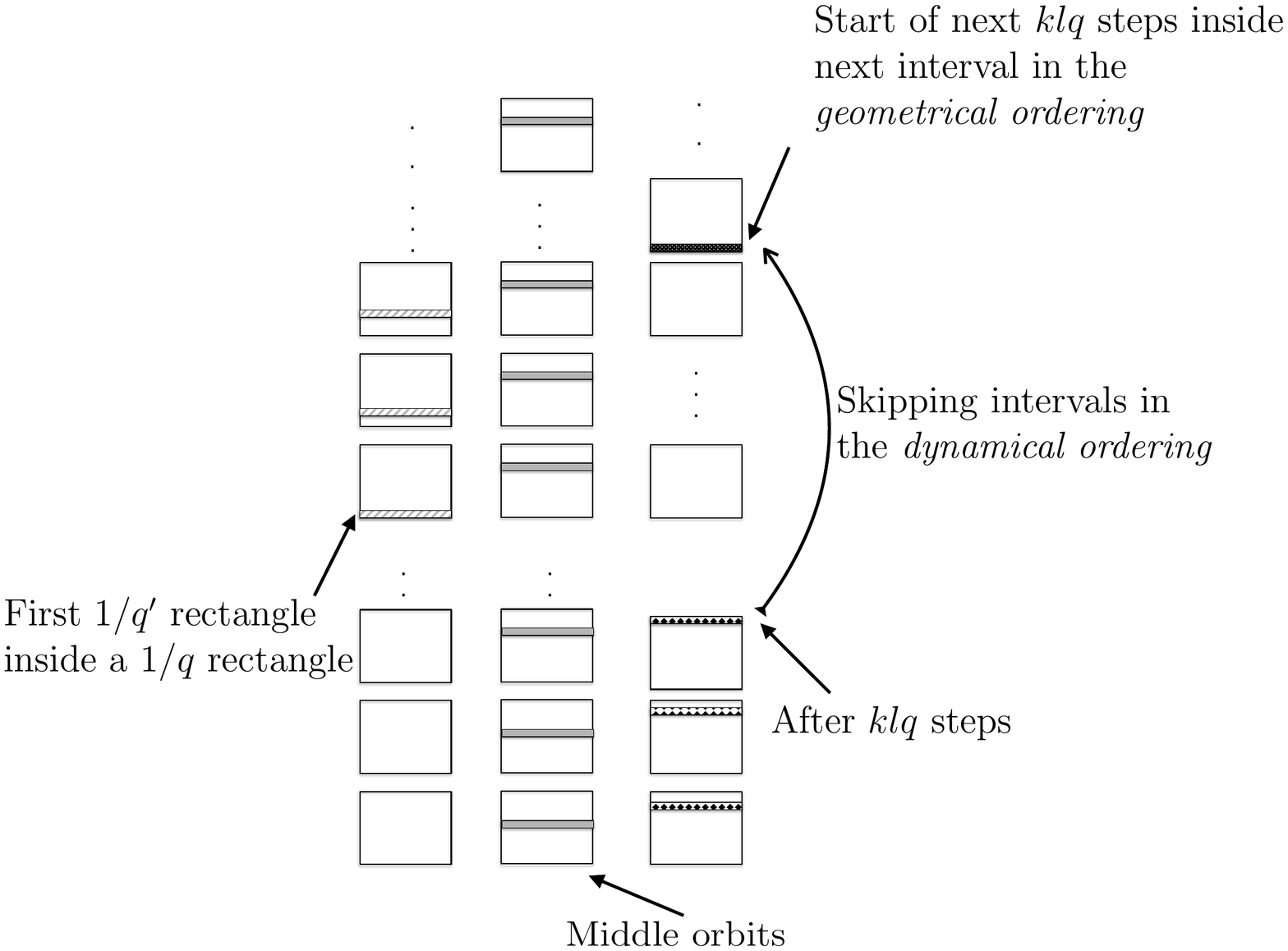}
\caption{{{\bf With vertical orientation}: Each column is the whole $1/q$ partition in dynamical ordering. The first shows the beginning trajectory of the $1/q'$ partition, the second shows one pass of many through the middle section, and the last column is the ending portion of the $1/q'$ trajectory. After the last step the end portion will jump some $1/q$ pieces in the dynamical ordering to begin a new trajectory.}}
\label{vertical list}
\end{figure}

{In the \emph{end} range (those $n$ for which $mkln+q-j_m+(kl-1)q\le n<(m+1)klq)$ the $\mcr^n_\beta J$  transits the 
intervals in places $0$ up to $j_m-1$ in the dynamical ordering on the $I^q_i$'s, ending up in the geometrically last 
subinterval of $I^q_m$. Finally we have that $\mcr^{(m+1)klq}_\beta J$ is the geometrically first subinterval of 
$I^q_{m+1}$. 

We illustrate this in figure \ref{vertical list} where we have three vertical columns. Each representing all of  the intervals of $\mci_q$ in their $\mcr_\alpha$-dynamical ordering. The transformation $\mcr_\alpha $ moves the larger rectangles up the columns, 
and $\mcr_\beta$ moves the flatter rectangles. The left hand column shows the beginning $\mci_{q'}$ intervals (starting with the first $1/q'$-interval for simplicity), the middle column shows one pass (of many) through the middle section and the last column shows the end portion. Note that the end portion ends in the rectangle just below where the beginning portion started. The black rectangle indicates the beginning of the next sequence after jumping over some of the $\mci_q$ rectangles in the geometric ordering. }


\pagebreak
\pagebreak

\subsubsection{With the partitions of $\zoo$}\label{on the interval}
We continue our examination of the path of $J=[0, 1/{q'})$ through the unit interval by considering it in light of the partition $\mci_{kq}$. In the discussion in this section we denote $I^q_i$ simply as $I_i$.

If $J'=J^{q'}_j$ is a subinterval of $I_i$   an application of $\mcr_\beta^{lq}$ moves $J'$ to the right ${lq\over klq^2}={1\over kq}$. In particular it moves it from being a subinterval of an element of $\mci_{kq}$ to the subinterval of the element of $\mci_{kq}$ adjacent on the right.

If $J'$ is the leftmost subinterval of $I_i$, then an application of $\mcr^{mlq}_\beta$ moves $J'$ over $m$
 elements of $\mci_{kq}$. For $0\le m<k, \mcr_\beta^{mlq}$ keeps $J'$ a subinterval of $I_i$, but 
 $\mcr_\beta^{klq}J'$ is the leftmost subinterval of $I_{i+1}$.

 We divide the atoms of $\mci_{kq}$ into $k$ ordered sets $w_0, \dots w_{k-1}$ where:
\begin{equation}
w_j=\la  I^{kq}_{j+tk}:0\le t<q\ra.\label{first w's}
\end{equation}
Thus each $w_j$ is the orbit under $\mcr_\alpha$ of $\hoo{j}{kq}$. We can view $w_j$ as a word of length $q$ in the alphabet $\mathcal I_{kq}$.


 We now follow our original interval $J$ through the $w_j$ under the iterates of $\mcr_\beta$. The $\mcr_\beta$--orbit
 of $J$ has length $q'$. The first $klq$ iterates are straightforward. Any number $n$ less than $klq$ can be 
 written  in the form $mlq+sq+t$ where $0\le m<k$, $0\le s<l$ and $0\le t<q$.  By the remarks in the previous 
 paragraph, $\mcr^n_\beta J$ is a subinterval of the $t^{th}$ element $w_m$. 
 We could then write the $\mci_{kq}$-name of any point in $J$ in the first $klq$ iterates as 
\begin{equation}\label{first time for everything}
w_0^lw_1^lw_2^l\dots w_{k-1}^l
\end{equation}

Applying $\mcr^{klq}_\beta$ to $J$ makes it be the geometrically first subinterval of $I_1$. The pattern above  would repeat itself with respect to the partition 
$\mci_{kq}$, were it not for the fact that the $w_j$'s start in $I_0$, and hence $\mcr_\beta^{klq}J$ is a subinterval of the $j_1^{st}$ element of $w_0$. Thus we must use $q-j_1$ applications of $\mcr_\beta$ to bring $\mcr_\beta^{klq}J$ back to $I_0$. Then $(l-1)q$ more applications carries $\mcr_\beta^{klq+(q-j_1)}J$ through $l-1$ copies of $w_0$, and $j_1$ additional applications brings it back to $I_1$ as a subinterval of an interval in the middle of $w_1$.  

In summary,  using $q-j_1$ iterates we are back in $I_0$, but in the first interval in $w_1$. Using $(l-1)q$ more iterates carries us through $l-1$ copies of $w_1$, and a $j_1$ more put us into the $j_1^{st}$ element  of $w_2$ and so on.
Hence we can write the $\mci_{kq}$ name of any point in $J$ in iterates 
$klq\le n<2klq$ in the form:
\begin{equation}
\label{first on and off ramps}
b_0^{q-j_1}w_0^{(l-1)}e_0^{j_1}b_1^{q-j_1}w_1^{l-1}e_1^{j_1}b^{q-j_1}w_2^{l-1}e_2^{j_1}\dots w_{k-1}^{l-1}e_{k-1}^{j_1}
\end{equation}
where $b_j^{q-j_1}$  and $e^{j_1}_j$  are the last $q-j_1$ elements of $w_j$  and the  first $j_1$-elements of $w_j$.

The interval $\mcr^{2klq}_\beta J$ is the leftmost subinterval of $I_2$ of length $1/q'$. This is a subinterval of the $j_2^{nd}$ element of the word $w_0$. Here the pattern repeats itself, but with $j_2$ playing the role of $j_1$ and so on.

Inductively one  shows that the $\mci_{kq}$-names of any element of $J$ have the following form between interates $mklq$ and $(m+1)klq$ with $1\le m<q$:
\begin{equation}\label{pass m}
b_0^{q-j_m}w_0^{l-1}e_0^{j_m}b_1^{q-j_m}w_1^{l-1}e_1^{j_m}b^{q-j_m}w_2^{l-1}e_2^{j_m}\dots w_{k-1}^{l-1}e_{k-1}^{j_m}
\end{equation} 
Since $j_0=0$, we can make equation \ref{first time for everything} a special case of equation \ref{pass m} by replacing the first copy of each $w_i$ by $b_i^{q-j_0}=b_i^{q}$. Using Remark \ref{noise}, we see that we get an isomorphic symbolic system.

Summarizing, any point in $J$ can have its $\mci_{kq}$-name written in the form 
\begin{equation}\label{first occurrence}
w=\prod_{i=0}^{q-1}\prod_{j=0}^{k-1}(b_j^{q-j_i}w_j^{l-1}e_j^{j_i}) 
\end{equation}
where the $b$'s and $e$'s have the definition given above. 

The letters occurring in particular locations of $w$ are in a specific one-to-one correspondence with intervals in the partition $\mci_{q'}$. So, for example, the second $b$ in the third string of the form $b^{q-j_3}$ uniquely labels an interval in $\mci_{q'}$ and so on. We will heavily use the correspondence between the letters of $w$ and the corresponding intervals. 

Another observation is that if we omit the subscripts on the $b$'s and $e$'s we can still decode this correspondence using the arithmetical properties of their exponents. 
Omitting the subscripts we see that the word $w$ in equation \ref{first occurrence} takes the form 
$\mcc(w_0, \dots w_{k-1})$.
\begin{remark}
By Lemma \ref{stabilization of names 1},  for all but $3/l$ portion of $x\in \zoo$, the $[-q,q]$-name of $x$ with respect to $\mci_{kq}$ under $\mcr_\beta$ is the same as the $[-q,q]$-name of $x$ with respect to $\mci_{kq}$ with under $\mcr_\alpha$.
\end{remark}


\subsubsection{What this means for the $\tau_n$}
The partition $\xi_{n}$ slices  $\mca$ into $s_n$ equal height horizontal strips  and into $q_n$ vertical strips over the atoms of $\mci_{q_n}$. We examine the analysis in the previous section taking $\alpha=\alpha_{n}$ and $\beta=\alpha_{n+1}$ (so $p=p_n, q=q_n, p'=p_{n+1}$ and $q'=q_{n+1}$). 


The action of $\rid{\alpha}$ on $\mca$ exactly mimics the action of $\mcr_{\alpha}$ on the $x$-axis. Hence on each horizontal  strip of 
$\mci_{k_nq_n}\otimes \mci_{s_{n+1}}$ of the form $\mci_{k_nq_n}\times \hoo{s}{s_{n+1}}$ we get the same analysis of the comparisons of $\rid{\alpha_{n+1}}$
 and $\rid{\alpha_n}$,  as we did in comparing $\mcr_{\alpha_n}$ and $\mcr_{\alpha_{n+1}}$.

The labeling of $\mci_{q_{n+1}}$ into words in equation \ref{first occurrence} can be copied over to $\mci_{q_{n+1}}\times \hoo{s}{s_{n+1}}$ and reflects the names of the   partition $\mci_{k_nq_n}\otimes \mci_{s_{n+1}}$ with respect to the action  of $\rid{\alpha_{n+1}}$. Doing this labels some of the atoms of $\mci_{q_{n+1}}\otimes \mci_{s_{n+1}}$ with some $b$'s and $e$'s, corresponding to the boundaries of the words of the type that occur in equation \ref{first occurrence}. Explicitly, if $a$ is an atom of $\mci_{q_{n+1}}$ labelled with a $b$ (respectively $e$) in equation \ref{first occurrence} then all of the atoms in $\{a\}\times \mci_{s_{n+1}}$ are labelled with a $b$ (respectively $e$).


\begin{definition}\label{definitions of B's and E's} Let $B_0=E_0=\emptyset$ and  let $B_{n+1}$ (respectively $E_{n+1}$) be the union of $B_n$ with the set of $x\in \mca$ that occur in
an atom of  $\mci_{q_{n+1}}\otimes \mci_{s_{n+1}}$  labelled with a $b$ (respectively an $e$). Let $B_{n+1}'=B_{n+1}\setminus B_n$ and $E'_{n+1}=E_{n+1}\setminus E_n$.
 \end{definition}
%
%
For $n>0$, measure of $B'_{n+1}\cup E'_{n+1}$ is $1/l_{n}$. 
 Moreover the collection of $x\in \mca$ whose $[-q,q]$-name with respect to the partition $\mci_{k_nq_n}\otimes \mci_{s_{n+1}}$ under $\rid{\alpha_n}$ is the same as its $[-q,q]$-name with under the action of $\rid{\alpha_{n+1}}$ has measure at least $1-3/l_n$. 

We now define \hypertarget{Gamma_n}{
\begin{equation}\label{def of gamma_n}\Gamma_n=\{x\in X:\mbox{for all }m>n, x\mbox{ does not occur in }Z_m(B'_m\cup E'_m)\}.
\end{equation}}
Then $\Gamma_{n+1}\supseteq \Gamma_{n}$. Since $\sum{1\over l_n}<\infty$,  the Borel-Cantelli Lemma implies that for almost every $x\in \mca$ there is an $m$ for all $n>m, x\in \Gamma_n$.

We note that the sets $B_n, E_n$ and $\Gamma_n$ correspond to the $B_n$, $E_n$ and $G_n$ in Section \ref{abstract reps of periodic}. The $B_n', E_n'$ are those $x$ that are first labelled $b$ or $e$ at stage $n$.
\subsection{The factor $\mck$ is a rotation of the circle}
As a warmup for the  symbolic representation of Anosov-Katok diffeomorphisms of surfaces, we are  in a position to give a one-dimensional representation of the circular system $\mck$ given in Definition \ref{first appearance of circle factor}. Let $\la k_n, l_n:n\in\nn\ra$ be a sequence of numbers such that $k_n\ge 2$ and $\sum 1/l_n<\infty$.  Let $p_n$ and $q_n$ be defined as in equations \ref{pn and qn}, $\alpha_n=p_n/q_n$. Then the sequence of $\alpha_n$ converge to an irrational $\alpha$. 

\begin{theorem}\label{mck}
Let $\nu$ be the unique non-atomic shift-invariant measure on $\mck$. Then 
\[(\mck, \mcb, \nu, sh)\cong(S^1, \mcd, \lambda, \mcr_\alpha)\]
where $\mcr_\alpha$ is the rotation of the circle $S^1$ and $\lambda$ is Lebesgue measure.
\end{theorem}
In the forthcoming \cite{global_structure} we give a completely different (short) algebraic proof of this result. We give a geometric proof here because it gives more information that is used in \cite{non-class}.

\pf Recall that the construction sequence for $\mck$ consists of $\mcw_n=\{w_n\}$ where $w_0=*$ and $w_{n+1}=\mcc(w_n, \dots w_n)$.

Let $X=[0,1)$ and take $\mu=\lambda$.
Define a sequence of uniform periodic processes that converge to $\mcr_\alpha$, by taking the $n^{th}$ periodic process $\sigma_n$ to be the cyclic permutation of $\mci_{q_n}$ given by the dynamical ordering. The base of $\sigma_n$ is $J_n=I^{q_n}_0$ and the levels are 
\[\la I_0, \mcr_{\alpha_n} I_0, \mcr_{\alpha_n}^2I_0 \dots \mcr_{\alpha_n}^{q-1}I_0\ra.\] 

The periodic approximation $\sigma_n$ can be realized pointwise by the transformation $\mcr_{\alpha_n}$. Since the $\mcr_{\alpha_n}$'s converge to $\mcr_\alpha$, and $\sum{q_n/q_{n+1}}<\infty$ Lemma \ref{convergent approximations} shows that the $\sigma_n$'s converge (as periodic processes) to $\mcr_{\alpha}$.

We now follow the method described in Section \ref{abstract reps of periodic}. Let $\la G_n:n\in\nn\ra$ be as defined there. Let $\mcq_0$ be the trivial partition of $G_0$. Labelling elements of $G_0$ with letter ``$*$", we inductively define $\la B_n, E_n:n\in\nn\ra$ and show  that if $x\in G_n$ is in the bottom level  of $\sigma_n$, then the $\{*,b,e\}$-name of $x$ is $w_n$.

We start with $\sigma_0$, the trivial action on the tower that has one level, the partition $\mci_1$.  The next periodic process $\sigma_1$ is a single cycle of length $q_1$.  We take the first level of $\sigma_1$ to consist of a subset of $B_1$ followed by $l_0-1$ levels which we view as 
 levels of $\sigma_0$ concatenated without spacers, followed by a single level contained in $B_1$, followed by $l_0-1$ levels without spacers, followed by a single level contained in $B_1$ and so on. There are $q_1=k_0l_0$ many levels total, of these $k_0$ are equally spaced subsets of $B_1$. The set $E_1=\emptyset$.

If we label the levels that are not subsets of $B_1$ with $*$'s and levels of $B_1$ with $b$'s we get a string that starts with a $b$, is followed by $l_0-1$ $*$'s, followed by a $b$ and so on, $k_1$ many times. Keeping in mind that $p_0=0$ and $q_0=1$ (so every $j_i=0$) this has the form
\[\prod_{i<q_0}\prod_{j<k_1}b^{q_0-j_i}*^{l_0-1}e^{j_i},\]
as desired.

We now describe the induction step.  Passing from $\sigma_n$ to $\sigma_{n+1}$ we assume that the words corresponding to elements of the bottom level of $\sigma_n$ 
are given by $w_n$. Following the analysis of Section \ref{on the interval}, we can partition the levels of $\sigma_{n+1}$ into contiguous segments that have names $w_n^{l_n-1}$ interspersed with $b$'s and $e$'s yielding the word $\mcc(w_n, w_n, \dots w_n)$. We let $B_n$ be the levels newly labelled with $b$'s and $E_n$ the levels newly labelled with $e$'s. By Lemma \ref{stabilization of names 1}, the measure of $B_n\cup E_n$ is $1/l$.

Let $S\subseteq \mck$ be as in Definition \ref{def of S}.   By Lemma \ref{dealing with S}, $\nu(S)=1$ and for all $s\in S$, there is an $N$ for all $n\ge N$ there are $a_n, b_n\ge 0$ such that $s\rest[-a_n, b_n)=w_n$.\footnote{The $s$'s for which $a_n$ and $b_n$ do not exist are those $s$ for which $s(0)$ is in the boundary portion of $w_m$ for some $m\ge n$.} For such an $s$ and $n$, let $r_n(s)=a_n$. We interpret $r_n(s)$ as the position of $s$'s ``$0$" in $w_n$. 

Supposing that $r_n(s)$ exists, we define 
{\begin{equation}\label{def of rhon}
{\rho}_n(s)={i\over q_n}
\end{equation}}
iff 
$I^{q_n}_i$ is the $r_n(s)^{th}$ interval in the dynamical ordering of $\mci_{q_n}$. (This  is equivalent to $i=j_{r_n(s)}$).\footnote{Thus $r_n$ and $\rho_n$ both have the same subset of $S$ as their domain and contain the same information. They map to different places $r_n:S\to \nn$, whereas $\rho_n:S\to\zoo$ and is the left endpoint of the $r_n^{th}$ interval in the dynamical ordering.} Equivalently, since 
the $r_n^{th}$ interval in the geometric ordering is $I^{q_n}_{p_nr_n(s)}$:
	\begin{equation*}
	i\equiv p_nr_n(s)\mod{q_n}
	\end{equation*}
Thus $\rho_n(s)$ is the left endpoint of the $r_n(s)^{th}$ interval in the periodic process $\sigma_n$.

Because the $r_{n+1}(s)^{th}$ letter in $w_{n+1}$ is in the  $r_n(s)^{th}$ position in a copy of $w_n$, we see that the 
$r_{n+1}(s)^{th}$ interval in the dynamical ordering of $\mci_{q_{n+1}}$ is a subinterval of the $r_n(s)^{th}$ interval in the dynamical ordering of $\mci_{q_n}$. It follows that 
\[\rho_{n+1}(s)\ge \rho_n(s)\]
and that 
\[|\rho_{n+1}(s)-\rho_n(s)|<1/q_n.\]
Since $\sum_n 1/q_n<\infty$, the sequence $\la\rho_n(s):n\in\nn\ra$ is Cauchy. We define
\[\phi_0(s)=\lim_n\rho_n(s).\]
It is easy to check that $\phi_0(sh(s))=\mcr_\alpha(\phi_0(s))$, and hence by the unique ergodicity of the measure $\nu$ on $S$
\[(\mck, \mcc, \nu, sh)\cong([0,1), \mcb, \lambda, \mcr_\alpha).\]
This finishes the proof.\qed

The following is immediate from the proof of Theorem \ref{mck} :

\begin{prop}\label{Dn and rn}
For $x\in [0,1)$ let  $D_n(x)=j$ if $x$ belongs to the $j^{th}$ interval in the dynamical ordering of $\mci_{q_n}$ (or equivalently, $D_n(x)=j$ if  $x\in I^{q_n}_{jp_n}$). Then for all $s\in S$ and all large enough $n$:
\[r_n(s)=D_n(\phi_0(s)).\]
\end{prop}

\subsection{A symbolic representation of the abstract Anosov-Katok systems}\label{a symb rep for AK}


We now give a symbolic representation of the transformations built by the version of the Anosov-Katok technique as described in Section \ref{AK method}. Our symbolic representation will consist of the names of points in $\mca$ with respect to a generating partition $\mcq$ that we build in section \ref{abstract reps of periodic}, with the addition of a systematic method of assigning $b$'s and $e$'s.  To find the names we compute them with respect   to the periodic processes $\tau_n$ and show that for  every $k$ and almost every point $x$, the $[-k,k]$-name of $x$ with respect to  $\mcq$ and $\tau_n$ stabilizes for large $n$.


If $\mcq$ is a partition that is refined by the levels of the towers of a periodic process $\tau$, then the $\mcq$-names of any pointwise realization of $\tau$ are constant on the levels of the tower. Hence we can view these names as naming the levels themselves in the periodic orbits of the action of $\tau$ on    various towers. We call the resulting collection of names the $(\tau, \mcq)$-names.

We begin our exercise by fixing an arbitrary partition $\mcq^*$ of $X$ that is refined by the partition $\zeta_n$ and comparing the $\mcq^*$-names of points under $\tau_n$ and $\tau_{n+1}$. 

\begin{remark}
For each $n$ we take the base of the $s^{th}$-tower in  the periodic process $\tau_n$ to be $Z_n(R^n_{0,s})$. It will follow that the word giving the 
$\mcq^*$-names for the $s^{th}$-tower is the same as the word consisting of the $(\rid{\alpha_n}, Z_n^{-1}(Q^*)$-names of the tower based at $R^n_{0,s}$ (which  is the first atom of $\xi_n$ lying on the $s^{th}$ horizontal strip of $\xi_n$).
\end{remark}

To compute the $(\tau_{n+1}, \mcq^*)$-names, we copy $\mcq^*$ to $\mca$ via $Z^{-1}_n$ to get a partition 
$\mcp=_{def}Z^{-1}_n\mcq^*$. Because 
\[\tau_{n+1}=Z_n(h_{n+1}\rid{\alpha_{n+1}}h_{n+1}^{-1})Z_n^{-1},\]
this reduces the problem of finding $(\tau_{n+1}, \mcq^*)$-names to that of
 computing  the $(h_{n+1}\rid{\alpha_{n+1}}h_{n+1}^{-1},\mcp)$-names of the towers whose levels constitute 
 the partition 
 $Z^{-1}_n\zeta_{n+1}$. Since $Z^{-1}_n\zeta_n=\xi_n$  and $h_{n+1}$ permutes the atoms of $\xi_{n+1}$, we see 
 that $Z_n^{-1}\zeta_{n+1}=\xi_{n+1}$. 
\bigskip

\hypertarget{moving rectangles}{{\bf\noindent How each rectangle moves:}} For notational simplicity, let $k=k_n, q=q_n$, $I_i=I^{q_n}_i$ and $J_j=I^{q_{n+1}}_j$.\footnote{We note that the behavior of $J_i$'s with respect to the partition $\mci_{k_nq_n}$ is  that of the transects we discussed in section \ref{on the interval}.} Fix a rectangle $R$ in $\xi_{n+1}$. We have two cases. The first case is that $h_{n+1}^{-1}R=R^{n+1}_{i,j}$ where $J_i$ is not the geometrically last $1/q_{n+1}$ subinterval of an interval  in 
$\mci_{kq}$. 

\begin{figure}[!h]
\centering
\includegraphics[height=.30\textheight]{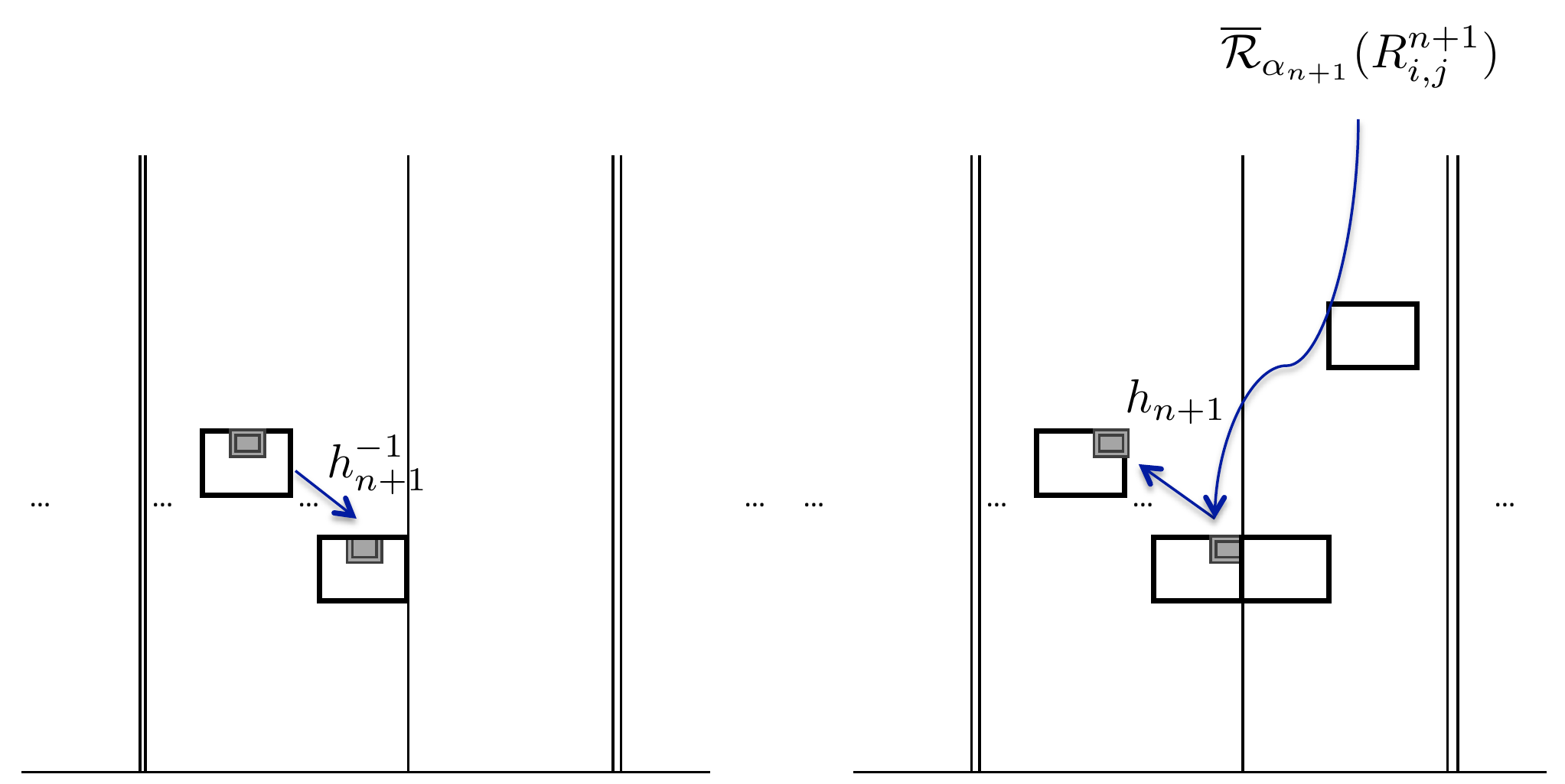}
\caption{Case 1}
\label{no jump rectangles}
\end{figure}

In this case know that $\rid{\alpha_{n+1}}$ and $\rid{\alpha_n}$ send $R^{n+1}_{i,j}$ to a rectangle whose base is a 
subinterval of the same element of $\mci_{kq}$.  Since $h_{n+1}$ commutes with $\rid{\alpha_n}$ and permutes the 
atoms of the  partition $\mci_{k_nq_n}\otimes \mci_{s_{n+1}}$ we see that $h_{n+1}\rid{\alpha_{n+1}}h_{n+1}^{-1}R$ is a 
subrectangle of the same member of $\mci_{kq}\otimes\mci_{s_{n+1}}$ as $\rid{\alpha_n}R$ is. In particular the $\mcp$-
name of $h_{n+1}\rid{\alpha_{n+1}}h_{n+1}^{-1}R$ is the same as the $\mcp$-name of $\rid{\alpha_n}R$.

The second case is when $J_i$ \emph{is} the geometrically last $1/q_{n+1}$ subinterval of an interval  in $\mci_{kq}$, then $\rid{\alpha_{n+1}}$ sends 
$R^{n+1}_{i,j}$
to the geometrically first subrectable of a new element $R'$ of $\mci_{kq}\otimes \mci_{s_{n+1}}$. 
Thus $h_{n+1}\rid{\alpha_{n+1}}h_{n+1}^{-1}R$ is a subset of $h_{n+1}(R')$. 
\begin{figure}[!h]
\centering
\includegraphics[height=.30\textheight]{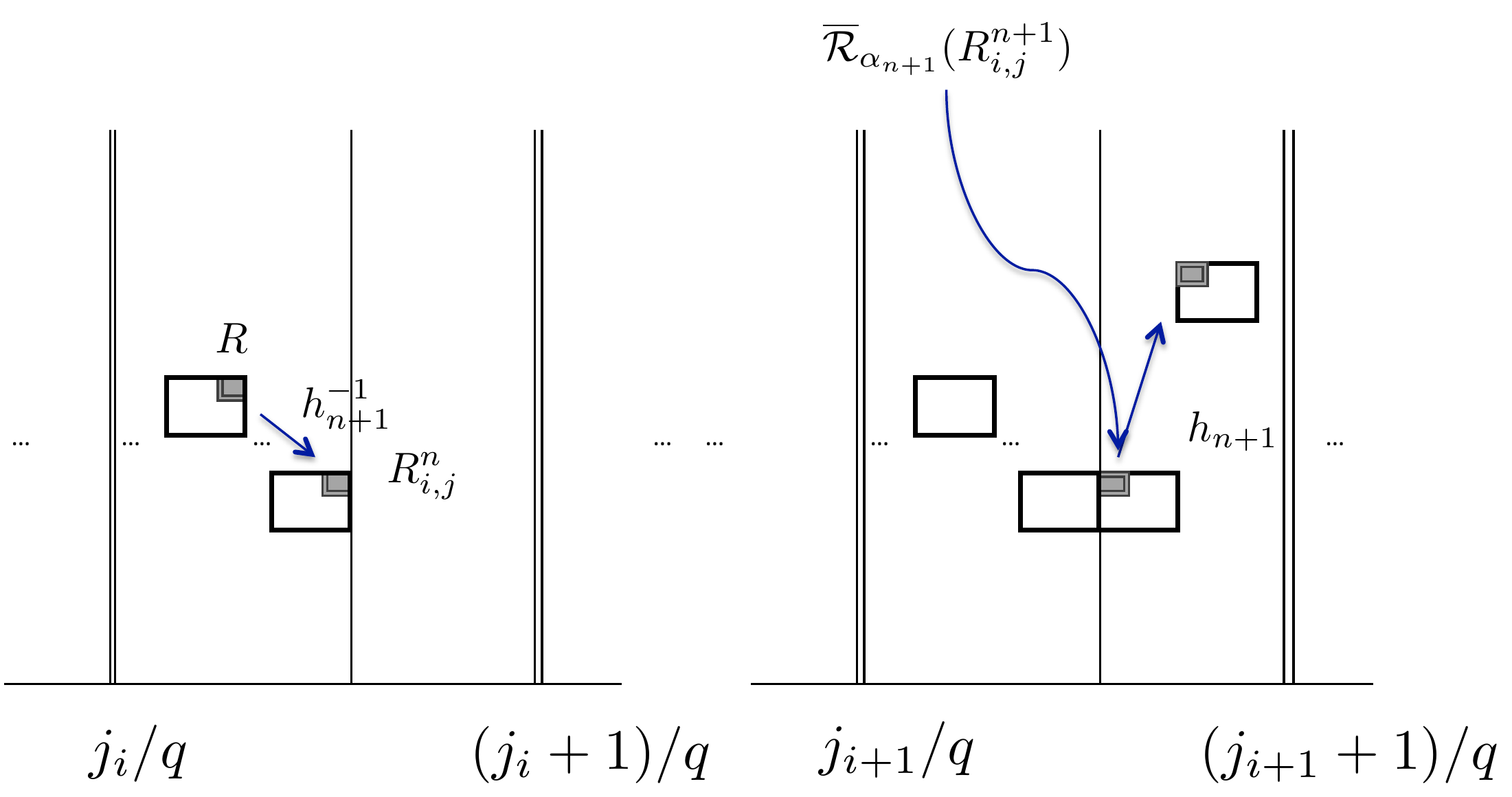}
\caption{Case 2}
\label{jumpy rectangles}
\end{figure}

\bigskip

{\bf \noindent How the tower moves:} Since the bases of the towers for $\tau_{n+1}$ are the sets  
$Z_{n+1}R^{n+1}_{0,s}$, the base for the $s^{th}$  tower for $h_{n+1}\rid{\alpha_{n+1}}h_{n+1}^{-1}$ is of the form 
\[Z_n^{-1}Z_{n+1}R^{n+1}_{0,s}=h_{n+1}R^{n+1}_{0,s}.\]
Computing:
\[(h_{n+1}\rid{\alpha_{n+1}}h_{n+1}^{-1})^t(h_{n+1}R_{0,s}^{n+1})=h_{n+1}\rid{\alpha_{n+1}}^tR^{n+1}_{0,s}.\]
Thus if $F_0, \dots F_{q_{n+1}-1}$ are the levels of the $s^{th}$ tower for $h_{n+1}\rid{\alpha_{n+1}}h_{n+1}^{-1}$, we see that $h_{n+1}^{-1}F_t=R^{n+1}_{i_t, s}$, 
 where $i_0=0$ and the sequence of intervals $\la J_{i_t}:t<q_{n+1}\ra$ is the orbit of $J_0$ under 
 $\mcr_{\alpha_{n+1}}$. 

In section \ref{on the interval}, we labelled the intervals $\la J_{i_t}:t<q_{n+1}\ra$ with $w$'s and $b$'s and $e$'s. From our discussion of how each \hyperlink{moving rectangles}{rectangle moves} we observe that for those $t$ where $J_{i_t}$ is  labelled with   a part of a $w$ the two transformations $h_{n+1}\rid{\alpha_{n+1}}h_{n+1}^{-1}$ and $\rid{\alpha_n}$ move $F_t$ to subrectangles of the same element of $\mci_{kq}\otimes\mci_{s_{n+1}}$ and hence the same element of $\mcp$.

Now, for $j<k, t<q$ and $s<s_{n+1}$, let $R_{j,t,s}$ be the rectangle 
\[\hoo{(j+tk)}{kq}\times [s/s_{n+1}, (s+1)/s_{n+1}).\]
This is the product of the $t^{th}$ interval in the $w_j$ of equation \ref{first w's} and the interval $\hoo{s}{s_{n+1}}$.
Let $u_{j,s}$ be the sequence of $\mcp$-names of the intervals
\begin{equation}\label{ujs}\la h_{n+1}(R_{j,t,s}):t<q\ra.
\end{equation} 
Then $u_{j,s}$ is the sequence of names  of the $h_{n+1}$-image of subrectangles of $\zoo\times \hoo{s}{s_{n+1}}$ 
taken along the transect whose horizontal intervals form the word $w_j$. Because $h_{n+1}$ commutes with $\rid{\alpha_n}$, we note the following:

\begin{remark} \label{coincidence} The word $u_{j,s}$ is the sequence of $\mcp$-names of the levels of the tower
 $\la {\rid{\alpha_n}^t}h_{n+1}(R^{n+1}_{i,s}):0\le t <q_n\ra$, for any $J_i\subseteq \hoo{j}{kq}$. 
\end{remark}

Following our analysis of the transects on $[0, 1)$, we can now describe the $\mcp$-name of the orbit of 
$F_0=h_{n+1}R^{n+1}_{0, s_0}$ under $h_{n+1}\rid{\alpha_{n+1}}h_{n+1}^{-1}$. The letter $t$  refers to the number 
of applications of $h_{n+1}\rid{\alpha_{n+1}}h_{n+1}^{-1}$ and thus the level of the tower.
\begin{enumerate}

\item To begin with there is a segment of $t$'s where the 
 $(h_{n+1}\rid{\alpha_{n+1}}h_{n+1}^{-1}, \mcp)$-name agrees with the $\rid{\alpha_n}$-name $u_{0,s}$. This segment has length  $lq_n$, as the intervals $J_{i_t}$ move to the right. At stage $t=lq_n$ these intervals cross a boundary for the $\mci_{kq}$ partition. At this point the name changes to $u_{1,s}$ repeated $l$-times. Then $u_{2,s}$ is repeated $l$ times and so on. This occurs $k$ times through the $u_{t,s}$ until $t=k_nl_nq_n-1$, where the $J_{i_t}$ becomes the geometrically last subinterval of $I_{q_n-1}$. Then $J_{i_{t+1}}$ is the geometrically first subinterval of $\hoo{1}{q_n}$.
\item We then have a segment where the transect is labelled with $q_n-j_1$ many $b$'s, after which $t=k_nl_nq_n +q_n-j_1$.

\item If $t=kq_n^l+ q_n-j_1$, then $J_{i_t}$ is a subinterval of the geometrically first interval in $\mci_{kq}$. In particular the name of $F_t$ is the first letter of $u_{0,s}$.
\item At this point the $h_{n+1}\rid{\alpha_{n+1}}h_{n+1}^{-1}$-names are the same as the $\rid{\alpha_n}$-names for $q(l-1)$ iterations yielding the name $u_{0,s}^{l-1}$. 
\item This is followed by a segment of length $j_1$ where the transect is labelled with $e$'s.
\item The pattern begins again with a portion of the tower where the transect is labelled with $b^{q-j_1}$, followed by $u_{1,s}^{l-1}$ followed by $e^{j_1}$. This is repeated for $u_{2,s}$, $u_{3,s}$ and so forth.
\item At stage $t=2k_nl_nq_n-1$,  $J_{i_t}$ is the geometrically last subinterval of $I_0$. This implies that
 $J_{i_{t+1}}$ is the geometrically first subinterval of $I_2$.
\item Here we get a block of $b$'s of length $q-j_2$ and the pattern described in items 3)-6) begins again with $j_2$ replacing $j_1$.

\item The pattern described in items 3)-8) repeat until we get to $t=klq^2-1$ at which point we have completed the period of $h_{n+1}\rid{\alpha_{n+1}}h^{-1}_{n+1}$.

\end{enumerate}

\begin{figure}[!h]
\centering
\includegraphics[height=.50\textheight]{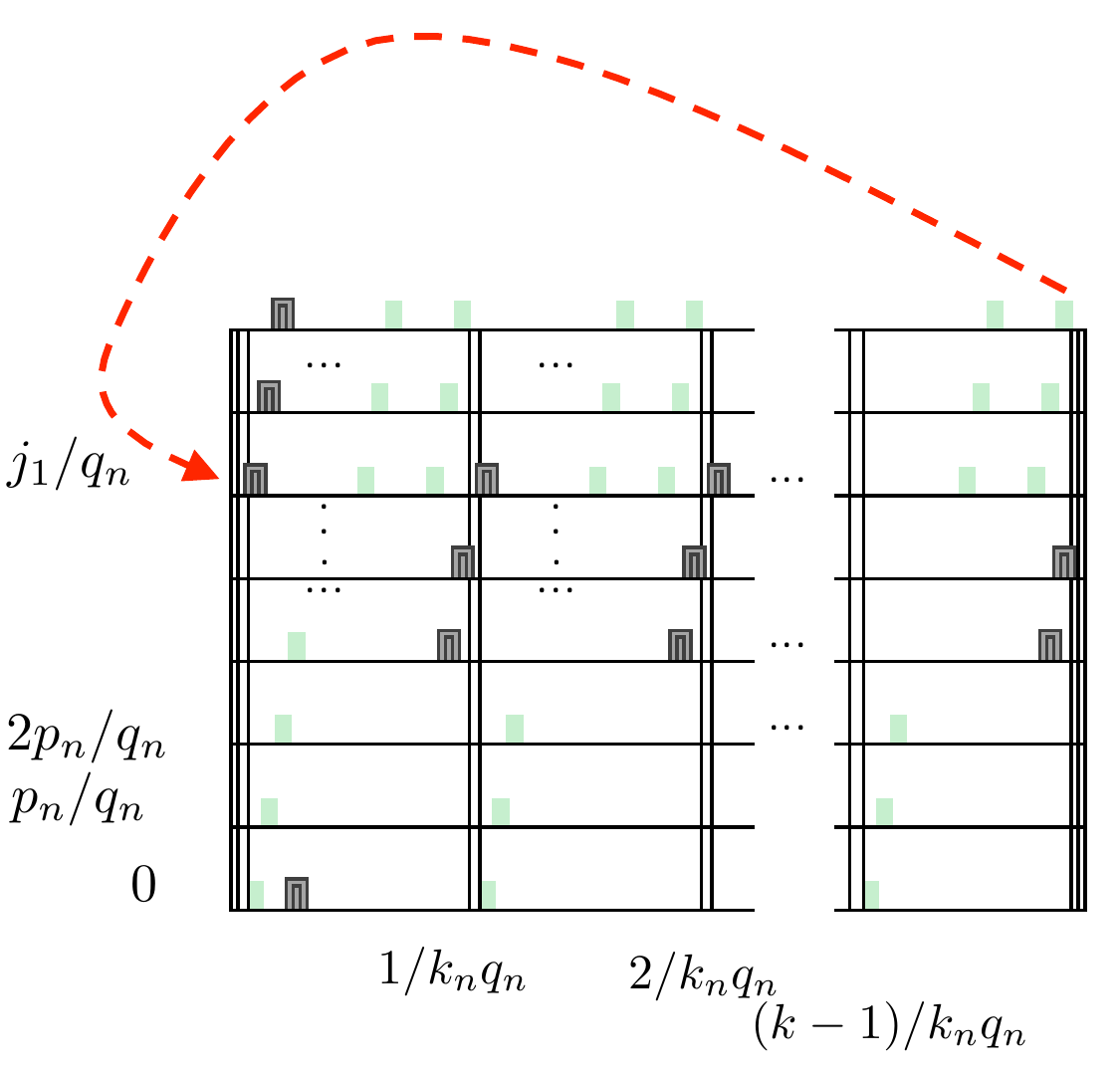}
\caption{A coarse diagram of the transects.}
\label{1/k}
\end{figure}

\begin{figure}[!h]
\centering
\includegraphics[height=.30\textheight]{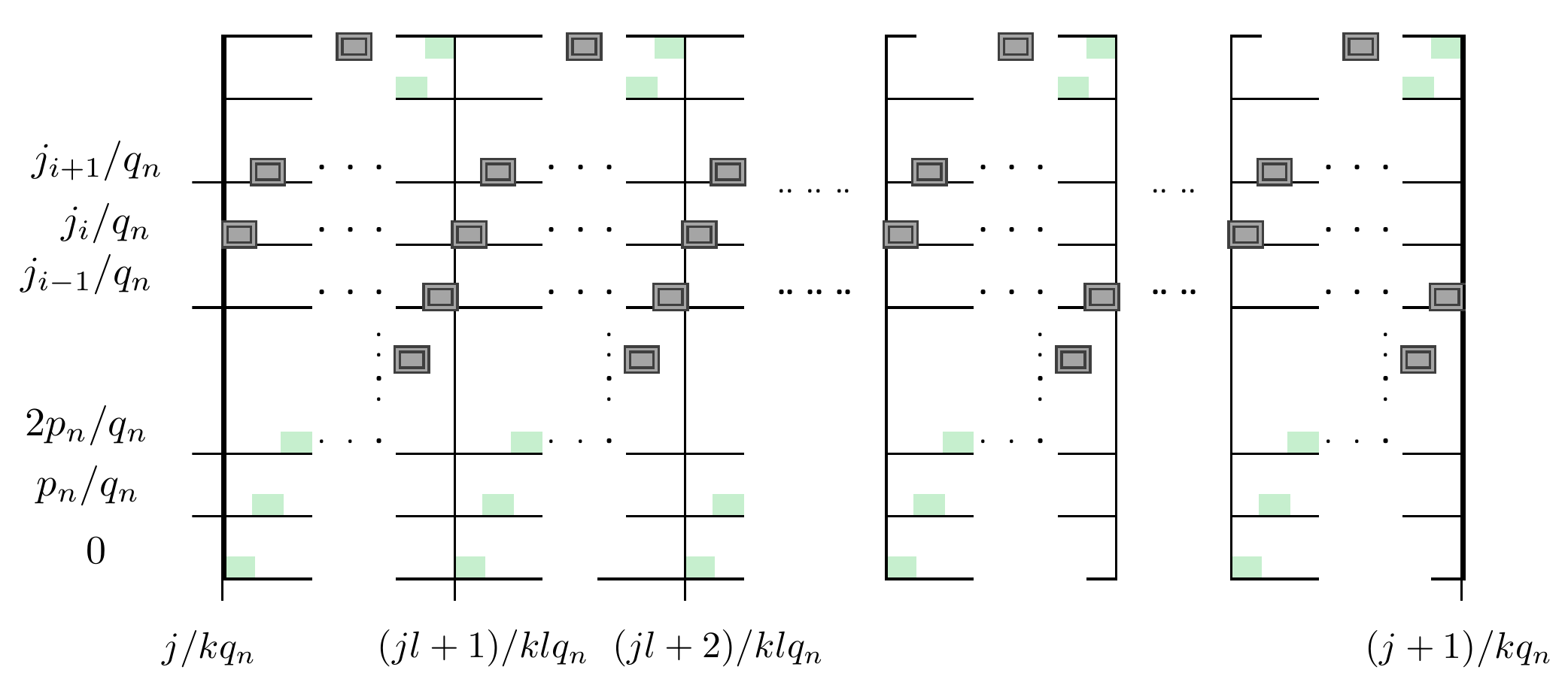}
\caption{A finer diagram of the transects. }
\label{1/l}
\end{figure}

We illustrate this with two diagrams. The levels of both figures  \ref{1/k} and \ref{1/l} are the rectangles of the partitions in the dynamical ordering with the $\mcr_{\alpha_n}$ and $\mcr_{\alpha_{n+1}}$ moving in the vertical direction. In figure \ref{1/k} we show how the $\xi_{n+1}$ transects move through $\mca$  in the global fashion at scale $1/k_nq_n$--neglecting the $1/k_nl_nq_n$ portions. The small light rectangles show the initial pass of a $\xi_{n+1}$ atom, with the more darkly shaded rectangles showing a later pass. 

In figure \ref{1/l} we magnify the first diagram to show the features at the $1/k_nl_nq_n$ scale. This is part of the darker rectangle transect from figure \ref{1/k} as it passes through a portion of width $1/k_nq_n$, going from $j/kq_n$ to $(j+1)/kq_n$ in $1/k_nl_nq_n^2$ increments.

We have shown:
\begin{theorem}\label{name computation}
Let $F_0$ be the base of a tower $\mct$ for $h_{n+1}\rid{\alpha_{n+1}}h_{n+1}^{-1}$. Suppose that $h_{n+1}^{-1}F_0=R_{0,s}$ Then the $\mcp$-names of $\mct$ agree with
\begin{equation*}u=_{def}\prod_{i=0}^{q-1}\prod_{j=0}^{k-1}(b^{q-j_i}u_{j,s}^{l-1}e^{j_i})
\end{equation*}
on the interior of $u$. 
\end{theorem}
From this we immediately get:
\begin{corollary}\label{substantial agreememnt}
Suppose that $J_i\subseteq [j/kq, (j+1)/kq-(1/q_{n+1}))$ with $j<kq$ and $R=R^n_{i,s}$. Then  the levels $\la\rid{\alpha_n}^th_{n+1}(R):t<q\ra$ coincide with the levels $\la(h_{n+1}\rid{\alpha_{n+1}}h_{n+1}^{-1})^th_{n+1}R:t<q\ra$ in the tower for $\tau_n$. In particular their $\mcq^*$-names agree.
\end{corollary}
\pf $J_i\subseteq [j/kq, (j+1)/kq-1/q_{n+1})$ is equivalent to $J_i$ being labelled with the first letter of a $w_m$ in a transect.\qed

In light of Remark \ref{stabilization of names 1}, we know:
\begin{corollary}\label{stabilization of names 2}
For a set of $x\in X$ having measure at least $1-3/l_n$, the $(\tau_{n+1},\mcq^*)$ and $(\tau_n,\mcq^*)$  names of $x$ agree on the interval $[-q,q]$.

\end{corollary}

We draw attention to the connection with  the sets \hyperlink{Gamma_n}{$\Gamma_n$} from equation \ref{def of gamma_n} (after Definition \ref{definitions of B's and E's}). Recall, that these are the collections of 
points that do not get labelled with $b$'s or $e$'s {for the first time} at some stage $m>n$. With this in mind the 
following corollary is clear:
\begin{corollary}\label{stable on gamma_n}
Suppose that $x\in \Gamma_n$ and $x$ is on level $t_n$ of a $\tau_n$-tower and the level $t_{n+1}$  of 
a $\tau_{n+1}$-tower. Let $w_n$ be the $\mcq^*$-name of $x$ with respect to $\tau_n$ and $w_{n+1}$ be the 
$\mcq^*$-name of $x$ with respect to $\tau_{n+1}$.  Then $w_{n+1}\rest[t_{n+1}-t_n, t_{n+1}+q_n-t_n)=w_n$.
\end{corollary}
This corollary is saying that for $x\in \Gamma_n$, the $\mcq^*$-name has stabilized on the interval of length $q_n$ corresponding to $x$'s position in a $\tau_n$-tower.
\bigskip

{\bf\noindent The symbolic representation:}
Let \[B=\{x\in X:\mbox{ for some }m\le n, x\in\Gamma_n\mbox{ and } Z^{-1}_mx\in B_m\},\] 
and 
\[E=\{x\in X:\mbox{ for some }m\le n, x\in\Gamma_n\mbox{ and }Z_m^{-1}x\in E_m\},\]
where the $B_n$'s and $E_n$'s are given by Definition \ref{definitions of B's and E's} and \hyperlink{Gamma_n}{$\Gamma_n$} is defined in equation \ref{def of gamma_n}.

Let $\{A_i:i<s_0\}$ be the partition $\zeta_0\rest(X\setminus B\cup E)$\footnote{Explicitly $\zeta_0\rest(X\setminus B\cup E)$ is given by the sets 
$Z[\zoo\times [s/s_0, (s+1)/s_0)]\setminus (B\cup E)$.} and 
\begin{equation}\label{def of Q}
\mcq=\{A_i:i<s_o\}\cup \{B, E\}.
\end{equation}

We:
\begin{enumerate} \item 
Compute the doubly infinite names of a typical point with respect to $\mcq$,
\item Show that ``\hyperlink{requirements 1-3}{{\bf Requirements 1-3"}} (occurring just before Lemma \ref{building hn from words}) imply $\mcq$  is a generator for the transformation $T=\lim \tau_n$, and
\item Show that the function sending an $x\in X$ to its $\mcq$-name has range in the set $S$  in Definition \ref{def of S}.
\end{enumerate}

The names will be in  the alphabet $\Sigma\cup \{b, e\}$ where $\Sigma=\{a_i:i<s_0\}$. Naturally a point $x\in X$ will get name $f\in (\Sigma\cup \{b,e\})^\poZ$ with $f(n)$ being  $a_i$ if $T^nx\in A_i$ and $f(n)$ being $b$ or $e$ if $T^nx\in B$ or $T^nx\in E$, respectively.


To give a complete description of the $(T,\mcq)$-names of points in $\bigcup \Gamma_n$ we go by induction on $n$. If $n=0$, then $\tau_0=id$, and $\Gamma_0=X\setminus(B\cup E)$. The $\tau_0$-names of $x\in \Gamma_0$ are simply the elements of $\Sigma$.

 Suppose that 
$F_0$ is the base of a  tower $\mct$ for $\tau_{n+1}$ and $Z_{n+1}R_0=F_0$ where $R_0=R^{n+1}_{0,s^*}$ for some $s^*<s_{n+1}$. Inductively assume that the $(\tau_n,\mcq)$-names of the towers with bases $Z_nR^{n}_{0,s}$ (with $s<s_n$) are $u_0, \dots u_{s_{n}-1}$. 
\begin{definition} Define a sequence of words $w_0, \dots w_{k_n-1}$ by setting  $w_j=u_s$ where 
	\[h_{n+1}(\hoo{j}{kq_n}\times \hoo{s^*}{s_{n+1}})\subseteq R^{n}_{0,s}.\] 
We will say that $(w_0, \dots w_{k_n-1})$ is the  sequence of $n$-words associated with $\mct$. 
\end{definition}

We define a circular system by inductively specifying the sequence $\la \mcw_n:n\in \nn\ra$. Let 
$\mcw_0=\{a_i:i<s_0\}$. 

Suppose that we have defined $\mcw_n$. Define $\mcw_{n+1}=\{\mcc_{n+1}(w_0, \dots w_{k_n-1}):$ $(w_0, \dots w_{k_n-1})$ is associated with a tower $\mct$  in $\tau_{n+1}\}$.  We will call $\la\mcw_n:n\in\nn\ra$ the {construction sequence} associated with the Ansov-Katok construction.  The following is worth noting.
\begin{prop}
Assume that the Anosov-Katok construction satisfies Requirements 1)-3) in section \ref{descs}. Then $\la\mcw_n:n\in\nn\ra$ is strongly uniform.
\end{prop}

\begin{theorem}\label{symbolic factor}
Suppose that $T$ is a limit of a sequence of Anasov-Katok periodic processes satisfying ``Requirements 1)-3)" with $ l_n$  growing fast enough. Then  almost all $x\in X$ have $\mcq$-names in 
 $\bk$, the circular system with construction sequence $\la\mcw_n:n\in\nn\ra$. In particular there is a measure $\nu$ on $\bk$ that makes $(\bk, \mcb, \nu, sh)$ isomorphic to the factor of $X$ generated by $\mcq$.
\end{theorem}

\pf Let $M$ be a positive integer. Then for almost all $x$ there is an $N$ so large that for all $n>N$, $x$ does not occur in the first or last $M$ levels of any tower in $\tau_n$.

Fix an arbitrary point $x\in \bigcup \Gamma_n$. Fix $n_0$ with $x\in \Gamma_{n_0}$ and consider $n\ge n_0$. Then 
$x$ belongs to a level of a tower $\mct$ of $\tau_{n}$. Without loss of generality we can assume that $x$ does not 
occur in the first or last $M$ levels of $\mct$. By Corollary \ref{stable on gamma_n}, if $x$ is in the $t_{n}^{th}$ level 
of a $\tau_n$-tower, then the $T$-name of $x$ agrees with the $\tau_n$-name of $x$ on the interval 
$[-t_n, q_n-t_n)$.

Applying Theorem \ref{name computation}, with $\mcp=Z_n^{-1}\mcq$, we see that a tower  $\mct$ for $\tau_{n}$ gets name
\begin{equation*}
w=\prod_{i=0}^{q_{n-1}-1}\prod_{j=0}^{k_{n-1}-1}(b^{q-j_i}w_{j}^{l-1}e^{j_i})
\end{equation*}
where $(w_0, \dots w_{k_{n-1}-1})$ is the sequence of words associated with $\mct$.
If $x$ is at  level $t_{n}$ then the $\mcq$-name of $x$ on the interval 
$[-t_{n}, q_{n}-t_{n})$ is 
\[\prod_{i=0}^{q-1}\prod_{j=0}^{k-1}(b^{q-j_i}w_{j}^{l-1}e^{j_i})=\mcc_{n}(w_0, \dots w_{k_{n-1}-1}).\]
Since $M<\min(t_{n}, q_{n}-t_{n})$, $x\rest[-M,M]$ is a subword of some word in $\mcw_{n}$.

Thus, for a typical $x$, every finite subinterval of the $(T,\mcq)$-name for $x$  is a subword of some 
$\mcw_{n+1}$. It follows that the factor of $(X, \mcb, \lambda, T)$ corresponding to the partition $\mcq$ is a factor of the uniform circular system we defined from the sequence of $\mcw_n$'s.


Conversely, by Lemma \ref{unique ergodicity}, the uniform circular system  $\bk$ determined by $\la \mcw_n:n\in\nn\ra$ is characterized as the smallest shift 
invariant closed set intersecting every basic open set $\la w\ra$ in $(\Sigma\cup \{b,e\})^\poZ$ determined by  
 some $w\in\mcw_n$. However each $w\in \mcw_n$ is represented on 
$\Gamma_n\cap \mct$ for some $\mct$, hence each $\la w\ra$ has non-empty intersection with the set of 
words arising from $(T,\mcq)$-names.\qed

We will need the following lemma that follows from   the proof of Theorem \ref{symbolic factor}. 
\begin{lemma}\label{being in S}
Let $a_0, b_0\in \nn$. Then for almost all $x$ and  all large $n$,  $x\in \Gamma_n$ and $x$ does not occur in the first $a_0$ levels or last $b_0$ levels of a tower of $\tau_n$.  In particular for almost all $x\in X$ there are $a>a_0, b>b_0$ such that
  the $\mcq$-name of $x$ restricted to the interval $[-a, b)$ belongs to  $\mcw_n$.
\end{lemma}

From Lemma \ref{being in S} and Theorem \ref{symbolic factor} we conclude:
\begin{corollary}\label{win}
For almost all $x\in X$ the $\mcq$-name of $x$ is in \hyperlink{def of s}{$S$}. In particular if $\nu$ is the unique non-atomic shift-invariant measure on $\bk$, then the factor $Y$ of $(X, \mcb, \lambda, T)$ generated by $\mcq$ is isomorphic to
 $(\bk, \mcb, sh, \nu)$. In particular, there is a unique non-atomic, shift-invariant measure on $Y$.

\end{corollary}

\bigskip

{\bf\noindent Generation.} To illustrate the potential difficulty, suppose that  each tower $\mct$ of each $\tau_n$ is associated to the same sequence of $n-1$-words 
then for every $x, t_1, t_2\in \zoo$, the $(\tau_n, \mcq)$-names of $(x, t_1)$ and $(x, t_2)$ are the same.  This would imply that $\mcq$ generates a proper factor of $X$.

Hence if $h_{n+1}$ does not vary enough on each horizontal strip of the form $\hoo{s}{s_{n+1}}$ the partition $\mcq$ may not generate. Fortunately Requirements 1-3 (stated just before Lemma \ref{building hn from words}) are sufficient conditions to guarantee generation. 

\begin{lemma}\label{generation x}
Suppose that for all $n$ the map sending a tower $\mct$ for $\tau_n$ to the sequence of $\mcq$-names associated to $\mct$ is a one-to-one function. Then $\mcq$ generates the transformation $T$.
\end{lemma}

\pf  Without loss of generality we can take $X=\mca$ and $Z$ to be the identity map. Since the $\la Z_n\xi_n:n\in\nn\ra$ is a decreasing sequence of partitions that generate the measure algebra, and $\mu(\Gamma_n)$ increases to 1, the atoms of $\la Z_n(\xi_n)\rest\Gamma_n:n\in\nn\ra$ also $\sigma$-generate the measure algebra.  Thus it suffices to show that each member of a $Z_n\xi_n\rest \Gamma_n$ belongs to the smallest translation invariant 
$\sigma$-algebra $\mcb$ generated by $\{\join_{i=-N}^{N}T^i(\mcq\cup \{B, E\}):N\in\nn\}.$


 Each $P\in Z_n\xi_n$ is the $t^{th}$ level of some tower $\mct$ for $\tau_n$. Let $w\in(\Sigma\cup\{b, e\})^{q_n}$ be the $(\tau_n, \mcq)$-name of $\mct$. Then the $j^{th}$ letter of $w$ determines an $S_{i_j}\in (\{B\}\cup \{E\}\cup \{A_i:i<s_0\})$. Since the $(\tau_n,\mcq)$-name of $\mct$ is correct on $\Gamma_n$, 
\[P\cap \Gamma_n\subseteq\bigcap_{0\le j<q_n}T^{t-j}(S_{i_j})\cap \Gamma_n.\]

On the other hand,  since the map sending towers to names $w$  is one-to-one we see
\[P\cap \Gamma_n\supseteq\bigcap_{0\le j<q_n}T^{t-j}(S_{i_j})\cap \Gamma_n,\]
  which is what we needed to show.\qed

It remains to show that the hypothesis of Lemma \ref{generation x} hold.

\begin{lemma} Suppose that our sequence of $h_n$'s satisfy the requirements 1-3 in section \ref{abstract AK}. Then $\mcq$ generates the transformation $T$.
\end{lemma}
\pf  We use requirement 3 to show inductively that for all $n\ge 1$, if $\mct$ and $\mct'$ are two $\tau_n$ towers then the $\mcq$-names associated with $\mct$ and $\mct'$ are different. 

For $n=0$ this is trivial. Suppose that it is true for $n$, we show it for $n+1$. Let $\mct$ and $\mct'$ be two towers and assume that they have bases $R^{n+1}_{0,s}$ and $R^{n+1}_{0,s'}$. By requirement 3, if $(j_0, \dots j_{k_n-1})$ and $(j'_0, \dots j'_{k_n-1})$ are the $k_n$ tuples associated with $s$ and $s'$, then they are distinct. Let $w_t$ be the $\mcq$-name associated with the $n$-tower with base $R^n_{0,t}$. By induction the $w_t$'s are distinct.  By Theorem \ref{name computation}, the $\mcq$-name of $\mct$ is 
$\mcc(w_{j_0}, \dots w_{j_{k_n-1}})$ and the $\mcq$-name of $\mct'$ is $\mcc(w_{j'_0}, \dots w_{j'_{k_n-1}})$. Since $(j_0, \dots j_{k_n-1})$ and $(j'_0, \dots j'_{k_n-1})$ are different we know $\mcc(w_{j_0}, \dots w_{j_{k_n-1}})$ and  $\mcc(w_{j'_0}, \dots w_{j'_{k_n-1}})$ are different.

It now follows from Lemma \ref{generation x} that $\mcq$ generates. \qed
\bigskip

{\bf\noindent Ergodicity of the Anosov-Katok systems:}
We can now show that abstract Anosov-Katok systems are ergodic and isomorphic to uniform circular systems.

\begin{theorem}\label{verifying unique ergodicity} Suppose that $(X, \mcb, \mu, T)$ is built by the Anosov-Katok method using fast growing coefficients and $h_n$'s satisfying requirements 1)-3). Let $\mcq$ be the partition defined in equation \ref{def of Q}.  Then the 
$\mcq$-names describe a strongly uniform circular construction sequence $\la \mcw_n:n\in\nn\ra$. Let $\bk$ be the associated circular system and  $\phi:X\to \bk$ be the map sending each $x\in X$ to its $\mcq$-name. {Then $\phi$ is one to one on a set of $\mu$-measure one. Moreover, }
there is a unique non-atomic shift-invariant measure $\nu$ concentrating on the range of $\phi$, and this measure is ergodic. In particular, $(X, \mcb, \lambda, T)$ is isomorphic to $(\bk, \mcb, \nu, sh)$ and is thus ergodic.
\end{theorem}

\pf Since the $l_n$-sequence grows fast, we know that the sequence of $\mcw_n$'s form a uniform construction sequence and hence there is a unique shift-invariant non-atomic measure $\nu$ on the set $S\subseteq \bk$ given in definition \ref{def of S} and $\nu$ is ergodic.
By Lemma \ref{being in S}, the range of $\phi$ is a subset of  the set $S$. Hence the factor determined by $\phi$ is isomorphic to $(S,\mcb, \mu, sh)$. In particular this factor is ergodic. 

Since the sequence of $h_n$'s satisfy requirements 1)-3) the partition $\mcq$ generates $X$. Hence $\phi$ is an isomorphism.  \qed

\begin{corollary}\label{representation}
If $T$ is a diffeomorphism of the disk built using the Anosov-Katok method of conjugacy satisfying requirements 1)-3) in section \ref{descs} then $T$ is measure theoretically isomorphic to a strongly uniform circular system.
\end{corollary}

\subsection{Tying it all together}
\label{symbolic summary}
In sections \ref{AK method} and \ref{symbolic representations of AK}, we have described a class of area 
preserving diffeomorphisms of the disk, annulus or torus. We have shown that, subject to some 
requirements (Requirements 1-3, in section \ref{abstract AK}), these transformations are ergodic and have a symbolic 
presentation in a particular form, that of uniform circular systems. 

The Anosov-Katok systems are built recursively depending on some data: some sequences of numbers $k_n, l_n, s_n$. Having been given these numbers, the final bit of data needed to determine the system is a sequence of permutations $h_n$ of the partitions $\xi_n$.
These permutations can be viewed as labeling the horizontal strips of the partition $\xi_n$ with bases of the towers from the previous periodic process. 

The numerical sequences and the labeling completely determine a construction sequence 
that is built recursively using an operator $\mcc$. The resulting sequence is uniform and circular and thus carries a unique non-atomic measure.

In our applications we take a different tack. We will view the results of this section as showing that uniform circular systems satisfying some minimal requirements are isomorphic to $C^\infty$-measure preserving transformations on the disk, annulus or torus.  Here is a converse to Corollary \ref{representation}.


\begin{theorem}\label{circular as smooth}
Suppose that $\la k_n, l_n, s_n:n\in \nn\ra$ are sequences of natural numbers tending to infinity such that  the $l_n$ grow sufficiently fast, {the $s_n$ grow to infinity} and $s_n$ divides both $k_n$ and $s_{n+1}$. 

Let $\la  \mcw_n:n\in\nn\ra$ be a circular construction sequence in an alphabet $\Sigma\cup \{b, e\}$  such that:
\begin{enumerate}
\item $\mcw_0=\Sigma$ and  for $n\ge 1,|\mcw_{n+1}|=s_{n+1}$, 
\item (Strong Uniformity)
For each $w'\in\mcw_{n+1}$, and $w\in\mcw_n$, if $w'=\mcc(w_0,\dots w_{k_n-1})$, then there are $k_n/s_n$ many $j$ with $w=w_j$.
\end{enumerate}
Then
\begin{enumerate}[A.)]
\item  If $\bk$ is the associated symbolic shift then there is a unique non-atomic ergodic measure $\nu$ on $\bk$.  

\item There is a $C^\infty$-measure preserving transformation $T$ defined on the torus (resp. disk, annulus) such that the system 
$(\mca, \mcb, \lambda, T)$ is isomorphic to $(\bk, \mcb, \nu, sh)$.
\end{enumerate}
\end{theorem}
Before we begin the proof of the theorem, we note that we do not know a \emph{a priori}  formulas for a growth rate for the $l_n$ that is sufficient for the conclusion of the theorem; the growth rate is determined inductively as described in the comments at the end of Section \ref{AK method}.
\medskip

\pf  We show that Lemma \ref{building hn from words} allows us to inductively construct   a sequence  $\la h_n:n\in\nn\ra$ that yields $\la \mcw_n:n\in\nn\ra$ as its construction sequence. Suppose that we have defined $\la h_{n^*}:n^*\le n\ra$.  From the definition of \emph{circular construction sequence} (Definition \ref{ccs})  we can find  $P_{n+1}\subseteq (\mcw_n)^{k_n}$ such that $\mcw_{n+1}$ is the collection of $w'$ such that for some sequence $(w_0, \dots w_{k_n-1})\in P_{n+1}, w'=\mcc(w_0, \dots w_{k_n-1})$. Enumerate $P_{n+1}$   as $w'_0, \dots w'_{s_{n+1}-1}$. Now apply Lemma \ref{building hn from words} to get $h_{n+1}$ from $w'_0, \dots w'_{s_{n+1}-1}$.

 We claim that the $\la h_n:n\in\nn\ra$ satisfy Requirements 1)-3). Requirement 1 is just that the $s_n$ go to infinity.  Requirement 2 follows from item 2 and requirement 3 follows since the words in $P_{n+1}$ are distinct.
\qed

The smooth transformation $T$ build in Theorem 
\ref{circular as smooth} is determined by the collections of words $\la P^T_n:n\in\nn\ra$ and our particular description of the Anosov-Katok construction.  
The words in 
$P^T_n$ determine the maps $h_n$ in the Anosov-Katok construction. Recall at the \hyperlink{fe}{end of section \ref{smooth AK} } we chose a summable sequence $\la \varepsilon_n:n\in\nn\ra$ such that $\varepsilon_n/4>\sum_{m>n}\varepsilon_m$ and a metric $d^\infty$ that determined the $C^\infty$-topology.  The sequence of $\varepsilon_n$ give estimates for the smooth approximations $h_n^s$ to $h_n$ and the sequence $\la S_n:n\in\nn\ra$ converging to $T$. Equation \ref{key to cont} shows that $ d^\infty(S_n,S_{n+1})<\varepsilon_n/4$. From this we observe that the
sequence $\la P^T_n:n\le M\ra$ determines an $\varepsilon_{M}$ neighborhood in which $T$ must lie. 

Conversely, different choices of $P_n$ give quite distant $h_n$'s and hence distant $h_n^s$ in the $C^\infty$-norm. 
We record this for use in applications.


\begin{prop}\label{continuous realization} Suppose that $\la \mcu_n:n\in\nn\ra$ and  $\la \mcw_n:n\in\nn\ra$ are construction sequences for two circular systems and $M$ is such that $\la \mcu_n:n\le M\ra=\la \mcw_n:n\le M\ra$. If $S$ and $T$ are the smooth realizations of the circular systems using the Anosov-Katok method given in this paper,  then  the $d^\infty$-distance between $S$ and $T$ is less than 
$\varepsilon_M$.
\end{prop}
\pf Given the circular construction sequences we have associated sequences $\la k_n^U, l_n^U, h_n^U, s_n^U:n\in \nn\ra$ and 
$\la k_n^W, l_n^W, h_n^W, s_n^W:n\in \nn\ra$ determining approximation $\la S_n:n\in\nn\ra, \la T_n:n\in\nn\ra$ to diffeomorphisms $S,T$. From the hypothesis we see
$\la k_n^U, l_n^U, h_n^U, s_n^U:n\le M\ra=\la k_n^W, l_n^W, h_n^W, s_n^W:n\le M\ra$. Thus $S_M=T_M$. By Remark \ref{going to a metric} and equation \ref{key to cont} we see that 
\begin{eqnarray*}
d^\infty(S_m,S)&<&\varepsilon_M/2\\
d^\infty(T_M,T)&<&\varepsilon_M/2.
\end{eqnarray*}
 It follows that $d^\infty(S,T)<\varepsilon_M$.
\qed

We end with a remark that seems relevant to the classification of diffeomorphisms of the torus up to conjugacy by homeomorphisms.

\begin{remark}
%
%

In Theorem 3.3 of \cite{FK} it is shown that there are exactly three ergodic invariant measures on the disk $D^2$ with respect to the Anosov-Katok diffeomorphism constructed there. There is one concentrating on the fixed point in the center, one concentrating on the boundary (where $T$ is a rotation) and one that gives every open set positive measure (Lebesgue measure). 

Because we are working on the torus, the top and bottom lines of $\mca$ are identified, rather than having one of them collapsed to a point. Thus we get two ergodic invariant measures. One concentrates on the ``equator" of the torus--the horizontal line corresponding to the top and bottom of the annulus we base our construction on. $T$ restricted to this line is the rotation $\mcr_\alpha$ ($\alpha=\lim \alpha_n$). The second invariant measure is Lebesgue measure.

If $\nu$ is an invariant ergodic measure that gives every open set positive measure then 
$\nu$ is Lebesgue measure. We can also prove this consequence in a different way. If $T$ is an Anosov-Katok diffeomorphism and $\bk$ is the circular system isomorphic to $T$, then there is a unique non-atomic invariant measure on $\bk$ (Lemma \ref{dealing with S}). If 
\[\phi:\bk\to \mathbb T^2\]
is the isomorphism, then the range of $\phi$ is $T$ invariant and has a unique non-atomic invariant measure. By considering the sets $G_n$ defined in equation \ref{Gn} one can establish that if $\nu$ is a $T$ invariant measure giving positive measure to every open set, then $\nu$ gives positive measure to the range of $\phi$. It follows that if $\nu$ is ergodic, then $\nu$ is Lebesgue measure.
\end{remark}
\subsection{Two Projects}
Here are two projects that we believe are of interest. 
The first is to extend the symbolic representation given in this paper to other versions of the Anosov-Katok construction--in particular to the twisted case, or to the constructions of weakly mixing transformations in \cite{FS}.

Secondly, the fact that weakly mixing transformations can be realized by the Anosov-Katok method suggest the possibility that a comeager collection of transformations (with respect to the  weak topology) could be realized by a method similar to the Anosov-Katok method.

\bibliography{citations}
\bibliographystyle{amsplain}
\end{document}